\documentclass[11pt]{amsart}
\usepackage[english]{babel}
\usepackage{amsthm} 
\usepackage{amssymb}
\usepackage{amsmath}
\usepackage{hyperref}
\usepackage{amsfonts}
\usepackage{verbatim}
\usepackage{stmaryrd}
\usepackage{fancyhdr}
\usepackage{mathrsfs}
\usepackage{graphicx}
\usepackage{pdfsync}
\newtheorem{theorem}{Theorem}[section]

\newtheorem{corollary}[theorem]{Corollary}

\newtheorem{lemma}[theorem]{Lemma}
\newtheorem{proposition}[theorem]{Proposition}
\newtheorem{remark}{Remark}

\numberwithin{equation}{section}

\newcommand{\N}{\mathbb{N}}
\newcommand{\R}{\mathbb{R}}

\newcommand{\lv}{\left[}
\newcommand{\rv }{\right]}
\newcommand{\rr}{\mathbb{R}}

\newcommand{\CC}{\mathbb{C}}
\newcommand{\cc}{\mathbb{C}}
\newcommand{\eps}{\varepsilon}

\newcommand {\be}{\begin{equation}}
\newcommand {\ee}{\end{equation}}
\newcommand {\ba}{\begin{array}}
\newcommand {\ea}{\end{array}}

\usepackage{varioref}
 
\def\Rom#1{\uppercase\expandafter{\romannumeral #1}}

\def\polhk#1{\setbox0=\hbox{#1}{\ooalign{\hidewidth \lower1.5ex\hbox{`}\hidewidth\crcr\unhbox0}}}

\begin{document}

\title[Some remarks on degenerate hypoelliptic Ornstein-Uhlenbeck operators]{Some remarks on degenerate hypoelliptic Ornstein-Uhlenbeck operators}

\author{M. Ottobre, G.A. Pavliotis \& K. Pravda-Starov}

\address{\noindent \textsc{Michela Ottobre, Department of Mathematics, Heriot-Watt University, Edinburgh EH14 4AS, UK}} 
\email{m.ottobre@hw.ac.uk}
\address{\noindent \textsc{Grigorios Pavliotis, Department of Mathematics, Imperial College London, 180 Queen's Gate,  London SW7 2AZ, UK}} 
\email{g.pavliotis@imperial.ac.uk}
\address{\noindent \textsc{Karel Pravda-Starov, IRMAR, CNRS UMR 6625, Universit\'e de Rennes 1, Campus de Beaulieu, 263 avenue du G\'en\'eral Leclerc, CS 74205,
35042 Rennes cedex, France}}
\email{karel.pravda-starov@univ-rennes1.fr}


\begin{abstract}
We study degenerate hypoelliptic Ornstein-Uhlenbeck operators in $L^2$ spaces with respect to invariant measures. The purpose of this article is to show how 
recent results on general quadratic operators apply to the study of degenerate hypoelliptic Ornstein-Uhlenbeck operators. We first show that some known results about the spectral and subelliptic properties of Ornstein-Uhlenbeck operators may be directly recovered from the general analysis of quadratic operators with zero singular spaces. 
We also provide new resolvent estimates for hypoelliptic Ornstein-Uhlenbeck operators.
We show in particular that the spectrum of these non-selfadjoint operators may be very unstable under small perturbations and that their resolvents can blow-up in norm far away from their spectra. Furthermore, we establish sharp resolvent estimates in specific regions of the resolvent set which enable us to prove exponential return to equilibrium.
\end{abstract}

\maketitle

\section{Introduction}

\subsection{Ornstein-Uhlenbeck operators}
We study Ornstein-Uhlenbeck operators
\begin{equation}\label{jen0}
P=\frac{1}{2}\sum_{i,j=1}^{n}q_{i,j}\partial_{x_i,x_j}^2+\sum_{i,j=1}^nb_{i,j}x_j\partial_{x_i}=\frac{1}{2}\textrm{Tr}(Q\nabla_x^2)+\langle Bx,\nabla_x\rangle, \quad x \in \rr^n,
\end{equation}
where $Q=(q_{i,j})_{1 \leq i,j \leq n}$ and $B=(b_{i,j})_{1 \leq i,j \leq n}$ are real $n \times n$-matrices, with $Q$ symmetric positive semi-definite. In the following, we denote $\langle A,B \rangle$ and $|A|^2$ the scalar operators
\begin{equation}\label{not}
\langle A,B\rangle=\sum_{j=1}^nA_jB_j, \quad |A|^2=\langle A,A\rangle=\sum_{j=1}^nA_j^2,
\end{equation}
when $A=(A_1,...,A_n)$ and $B=(B_1,...,B_n)$ are vector-valued operators. Notice that $\langle A,B\rangle\neq\langle B,A\rangle$ in general, since e.g., 
$\langle \nabla_x,Bx\rangle=\langle Bx,\nabla_x\rangle +\textrm{Tr}(B).$

We consider degenerate hypoelliptic Ornstein-Uhlenbeck operators for which the symmetric matrix $Q$ is not positive definite, but only positive semi-definite. 
In the recent years, these degenerate operators have been studied extensively, e.g.~\cite{lanco,farkas,lunardi1,lanco2,Lorenzi,Metafune_al2002}. We recall from these works that 
the assumption of hypoellipticity is characterized by the following equivalent assertions:

\medskip
 
\begin{itemize}
\item[$(i)$] The Ornstein-Uhlenbeck operator $P$ is hypoelliptic.
\item[$(ii)$] The symmetric positive semi-definite matrices
\begin{equation}\label{pav0}
Q_t=\int_0^{t}e^{sB}Qe^{sB^T}ds,
\end{equation}
with $B^T$ the transpose matrix of $B$, are nonsingular for some (equivalently, for all) $t>0$, i.e. $\det Q_t>0$.
\item[$(iii)$] The Kalman rank condition holds: 
\begin{equation}\label{kal1}
\textrm{Rank}[B|Q^{\frac{1}{2}}]=n,
\end{equation} 
where 
$$[B|Q^{\frac{1}{2}}]=[Q^{\frac{1}{2}},BQ^{\frac{1}{2}},\dots, B^{n-1}Q^{\frac{1}{2}}],$$ 
is the $n\times n^2$ matrix obtained by writing consecutively the columns of the matrices $B^jQ^{\frac{1}{2}}$, with $Q^{\frac{1}{2}}$ the symmetric positive semi-definite matrix given by the square root of $Q$.
\item[$(iv)$] H\"ormander's condition holds:
$$\forall x \in \rr^n, \quad \textrm{Rank } \mathcal{L}(X_1,X_2,...,X_n,Y_0)(x)=n,$$
with 
$$Y_0=\langle Bx,\nabla_x\rangle, \quad X_i=\sum_{j=1}^nq_{i,j}\partial_{x_j}, \quad i=1,...,n$$ 
 \end{itemize}

\medskip

\noindent
When the Ornstein-Uhlenbeck operator is hypoelliptic, i.e. when one (equivalently, all) of the above conditions holds, the associated Markov semigroup $(T(t))_{t \geq 0}$ has the following explicit representation
\begin{equation}\label{gen}
(T(t)f)(x)=\frac{1}{(2\pi)^{\frac{n}{2}}\sqrt{\det Q_t}}\int_{\rr^n}e^{-\frac{1}{2}\langle Q_t^{-1}y,y\rangle}f(e^{tB}x-y)dy, \quad t>0.
\end{equation}
This formula is due to Kolmogorov~\cite{kolmo}. On the other hand, the existence of an invariant measure $\mu$ for the Markov semigroup $(T(t))_{t \geq 0}$, i.e., a probability measure on $\rr^n$ verifying 
$$\forall t \geq 0, \forall f \in C_b(\rr^n), \quad \int_{\rr^n}(T(t)f)(x)d\mu(x)=\int_{\rr^n}f(x)d\mu(x),$$
where $C_b(\rr^n)$ stands for the space of continuous and bounded functions on $\rr^n$, is known to be equivalent~\cite[Sec.~11.2.3]{DaPrato}
to the following localization of the spectrum of~$B$,
\begin{equation}\label{pav1}
\sigma(B) \subset \mathbb{C}_-=\{z \in \mathbb{C} : \textrm{Re }z<0\}.
\end{equation}
When this condition holds, the invariant measure is unique and is given by $d\mu(x)=\rho(x)dx$, where the density with respect to the Lebesgue measure is
\begin{equation}\label{pav2}
\rho(x)=\frac{1}{(2\pi)^{\frac{n}{2}}\sqrt{\det Q_{\infty}}}e^{-\frac{1}{2}\langle Q_{\infty}^{-1}x,x\rangle},
\end{equation}
with 
\begin{equation}\label{pav3}
Q_{\infty}=\int_0^{+\infty}e^{sB}Qe^{sB^T}ds.
\end{equation}
The spectral and subelliptic properties of hypoelliptic Ornstein-Uhlenbeck operators have been studied in detail recently~\cite{lanco,lunardi1,Metafune_al2002}. We recall in particular that the spectrum of these operators in $L^p_{\mu}=L^p(\rr^n,d\mu)$ spaces with respect to invariant measures for $1 \leq p <+\infty$, has been explicitly described by Metafune, Pallara and Priola~\cite{Metafune_al2002}. Furthermore, the work by Farkas and Lunardi~\cite{lunardi1} provides optimal embeddings for the domains of hypoelliptic Ornstein-Uhlenbeck operators on $L^2_{\mu}$ spaces with respect to invariant measures, whereas global $L^p=L^p(\rr^n,dx)$ estimates of the elliptic directions were established by Bramanti, Cupini, Lanconelli and Priola~\cite{lanco}
for $1<p<+\infty$, see below the estimates (\ref{ellip}).

The purpose of this article is to show how recent results on general quadratic operators allow us to provide a unified approach for studying these problems.
We first show that some of the results mentioned above for Ornstein-Uhlenbeck operators may be completely or partially recovered from the general analysis of quadratic operators developed  in the works~\cite{HPS, OPPS,duke,karel}, see Propositions~\ref{prop1},~\ref{rt2} and~\ref{subest} in the next section. We also study the spectral stability of these operators under small perturbations. These operators are in general non-normal. It is for instance always the case when these hypoelliptic operators are degenerate, i.e., when the symmetric matrix $Q$ is not positive definite. When an operator is non-normal, it is well-known that its resolvent can blow-up in norm in unbounded regions of the resolvent set very far away for the spectrum~\cite{trefethen,trefethen2}. As recalled in the following, this type of phenomena is linked to the very strong instability of the spectrum of the operator under small perturbations. In the present work, we show that such phenomena occur for all non-normal hypoelliptic Ornstein-Uhlenbeck operators on $L^2_{\mu}$ spaces with respect to invariant measures (Theorems~\ref{esti} and~\ref{esti9}). On the other hand, we show that the hypoellipticity of Ornstein-Uhlenbeck operators still allows us to establish sharp resolvent estimates in specific regions of the resolvent set, whose geometry directly depends on the loss of derivatives with respect to the elliptic case in global subelliptic estimates satisfied by these operators (Theorem~\ref{esti1}). These resolvent estimates are then very useful to control and get sharp bounds for the associated semigroups which allow to establish results of exponential return to equilibrium (Theorem~\ref{return}). Related results for degenerate hypoelliptic Fokker-Planck operators are then given in Propositions~\ref{rt277} and \ref{return77}, and Corollary~\ref{return776}.

\subsection{Setting of the analysis}

Let $P$ be the Ornstein-Uhlenbeck operator defined in (\ref{jen0}). When $P$ is hypoelliptic and admits an invariant measure, we may associate to the operator $P$ acting on $L^2_{\mu}=L^2(\rr^n,d\mu)$, the quadratic operator $\mathscr{L}$ acting on $L^2=L^2(\rr^n,dx)$,
\begin{equation}\label{pav4.5}
\mathscr{L}u=-\sqrt{\rho}P\big((\sqrt{\rho})^{-1}u\big)-\frac{1}{2}\textrm{Tr}(B)u.
\end{equation}
We notice that the localization of the spectrum (\ref{pav1}) implies that 
$$\textrm{Tr}(B)<0,$$
since $B \in M_n(\rr)$. 
Recalling the notation (\ref{not}), a direct computation (see (\ref{pav15}) in Section~\ref{pav4}) shows that
\begin{multline}\label{pav5}
\mathscr{L} =-\frac{1}{2}|Q^{\frac{1}{2}}\nabla_x|^2+\frac{1}{8}|Q^{\frac{1}{2}}Q_{\infty}^{-1}x|^2-\Big\langle\Big(\frac{1}{2}QQ_{\infty}^{-1}+B\Big)x,\nabla_x\Big\rangle\\
=\frac{1}{2}|Q^{\frac{1}{2}}D_x|^2+\frac{1}{8}|Q^{\frac{1}{2}}Q_{\infty}^{-1}x|^2-i\Big\langle\Big(\frac{1}{2}QQ_{\infty}^{-1}+B\Big)x,D_x\Big\rangle, 
\end{multline} 
with $D_x=i^{-1}\nabla_x$, where $Q_{\infty}$ is the symmetric positive definite matrix (\ref{pav3}).
The operator $\mathscr{L}$ may be considered as a pseudodifferential operator
\begin{equation}\label{3}
\mathscr{L}u=q^w(x,D_x)u(x) =\frac{1}{(2\pi)^n}\int_{\R^{2n}}{e^{i(x-y) \cdot \xi}q\Big(\frac{x+y}{2},\xi\Big)u(y)dyd\xi},
\end{equation}
defined by the Weyl quantization of the quadratic symbol
\begin{equation}\label{pav6}
q(x,\xi)=\frac{1}{2}|Q^{\frac{1}{2}}\xi|^2+\frac{1}{8}|Q^{\frac{1}{2}}Q_{\infty}^{-1}x|^2-i\Big\langle\Big(\frac{1}{2}QQ_{\infty}^{-1}+B\Big)x,\xi\Big\rangle, \quad (x,\xi) \in \rr^{2n}.
\end{equation}
This explicit computation is performed in (\ref{pav17}) (Section~\ref{weyl}) by noticing that the Weyl quantization of the quadratic symbol
$x^{\alpha} \xi^{\beta}$, with $(\alpha,\beta) \in \N^{2n}$, $|\alpha+\beta|=2$, is the differential operator
\begin{equation}\label{forsym}
(x^{\alpha} \xi^{\beta})^w=\textrm{Op}^w(x^{\alpha} \xi^{\beta})=\frac{x^{\alpha}D_x^{\beta}+D_x^{\beta} x^{\alpha}}{2}. 
\end{equation}
Notice that the real part of the symbol $\textrm{Re }q \geq 0$ is a non-negative quadratic form
since $Q$, $Q_{\infty}$ and $B \in M_n(\rr)$. We know from~\cite{mehler} (p.~425) that the maximal closed realization of the operator $\mathscr{L}$, i.e., the operator on $L^2$ with domain
\begin{equation}\label{dom1}
D(\mathscr{L})=\{u \in L^2 : \ \mathscr{L}u \in L^2\},
\end{equation}
coincides with the graph closure of its restriction to the Schwartz space
$$\mathscr{L} : \mathscr{S}(\rr^n) \rightarrow \mathscr{S}(\rr^n).$$
Classically, to any quadratic form defined on the phase space 
$$q : \rr_x^n \times \rr_{\xi}^n \rightarrow \mathbb{C}, \quad n \geq 1,$$
is associated a matrix $F \in M_{2n}(\CC)$ called its Hamilton map, or its fundamental matrix, which is defined as the unique matrix satisfying the identity
\begin{equation}\label{10}
\forall  (x,\xi) \in \R^{2n},\forall (y,\eta) \in \R^{2n}, \quad q((x,\xi);(y,\eta))=\sigma((x,\xi),F(y,\eta)), 
\end{equation}
with $q(\cdot;\cdot)$ the polarized form associated to the quadratic form $q$, where $\sigma$ stands for the standard symplectic form
\begin{equation}\label{11}
\sigma((x,\xi),(y,\eta))=\langle \xi, y \rangle -\langle x, \eta\rangle=\sum_{j=1}^n\xi_j y_j-x_j \eta_j,
\end{equation}
with $x=(x_1,...,x_n)$, $y=(y_1,....,y_n)$, $\xi=(\xi_1,...,\xi_n)$, $\eta=(\eta_1,...,\eta_n) \in \cc^n$. 
We check in (\ref{pav20}) (Section~\ref{weyl}) that the Hamilton map of the quadratic form (\ref{pav6}) is given by
\begin{equation}\label{F}
F=\lv\begin{array}{cc}
 -\frac{i}{4}(QQ_{\infty}^{-1}+2B)&\frac{1}{2}Q \\
 -\frac{1}{8}Q_{\infty}^{-1}QQ_{\infty}^{-1}  & \frac{i}{4}(QQ_{\infty}^{-1}+2B)^T       
\end{array}\rv.
\end{equation}

In~\cite{HPS}, the notion of singular space was introduced by pointing out the existence of a particular vector subspace of the phase space, which is intrinsically associated to a quadratic symbol $q(x,\xi)$ and defined as the following finite intersection of kernels
\begin{equation}\label{h1bis}
S=\Big( \bigcap_{j=0}^{2n-1}\textrm{Ker}
\big[\textrm{Re }F(\textrm{Im }F)^j \big]\Big)\cap \rr^{2n},
\end{equation}
where $\textrm{Re }F$ and $\textrm{Im }F$ stand respectively for the real and imaginary parts of the Hamilton map $F$ associated to the quadratic symbol $q$,
$$\textrm{Re }F=\frac{1}{2}(F+\overline{F}), \quad \textrm{Im }F=\frac{1}{2i}(F-\overline{F}).$$
The works~\cite{HPS,OPPS,karel} show that this particular vector subspace of the phase space plays a basic role in the understanding of the properties of the quadratic operator $q^w(x,D_x)$. 
According to (\ref{F}) and (\ref{h1bis}), the singular space of the quadratic form (\ref{pav6}) reads as
\begin{multline}\label{pav7}
S=\big\{(x,\xi) \in \rr^{2n} : \forall 0 \leq j \leq 2n-1,\ QQ_{\infty}^{-1}(QQ_{\infty}^{-1}+2B)^jx=0, \\ Q(Q_{\infty}^{-1}Q+2B^T)^j\xi=0\big\},
\end{multline}
since $Q$ and $Q_{\infty}$ are symmetric matrices.
The key point in this article is to prove in (\ref{pav50}) (Section~\ref{singul}) that this singular space is actually equal to zero 
\begin{equation}\label{xc1}
S=\{0\}.
\end{equation}
A complete proof of this key property of hypoelliptic Ornstein-Ulhenbeck operators with invariant measures is given in Section~\ref{ele}.

For now, the next section shows how the recent results~\cite{HPS,OPPS,karel}  about general quadratic operators with zero singular spaces apply to the study of degenerate hypoelliptic Ornstein-Uhlenbeck operators or degenerate hypoelliptic Fokker-Planck operators, and how they relate to the previous works on this subject~\cite{anton,lanco,lunardi1,Metafune_al2002}.

\section{Statement of Main Results}

\subsection{Smoothing effect}
The first proposition shows that the semigroup associated to a hypoelliptic Ornstein-Uhlenbeck operator with an invariant measure satisfies smoothing and decay properties:

\medskip

\begin{proposition}\label{prop1}
Let 
$$P=\frac{1}{2}\emph{\textrm{Tr}}(Q\nabla_x^2)+\langle Bx,\nabla_x\rangle,\quad x \in \rr^n,$$ 
be a hypoelliptic Ornstein-Uhlenbeck operator, which admits the invariant measure $d\mu(x)=\rho(x) dx$. Then, the global solution to the Cauchy problem
\begin{equation}\label{p2}
\left\lbrace\begin{array}{c}
\partial_tv=Pv,\\
v|_{t=0}=v_0 \in L_{\mu}^2,
\end{array}\right.
\end{equation}
satisfies the following property
$$\forall t >0,\forall \alpha, \beta \in \mathbb{N}^n, \quad x^{\alpha}\partial_{x}^{\beta}v(t)=x^{\alpha}\partial_{x}^{\beta}(e^{tP}v_0) \in L_{\mu}^2.$$
\end{proposition}

\medskip

\begin{proof}
By using the fact that the quadratic operator
$$\mathscr{L}=q^w(x,D_x)=\frac{1}{2}|Q^{\frac{1}{2}}D_x|^2+\frac{1}{8}|Q^{\frac{1}{2}}Q_{\infty}^{-1}x|^2-i\Big\langle\Big(\frac{1}{2}QQ_{\infty}^{-1}+B\Big)x,D_x\Big\rangle,$$
has a Weyl symbol with a zero singular space $S=\{0\}$ and a non-negative real part $\textrm{Re }q \geq 0$,
we can deduce from~\cite{HPS} (Theorem~1.2.1) that the evolution equation associated to the accretive operator $\mathscr{L}$, 
$$\left\lbrace\begin{array}{c}
\partial_tu+\mathscr{L}u=0,\\
u|_{t=0}=u_0 \in L^2(\rr^n,dx),
\end{array}\right.$$
is smoothing in the Schwartz space $\mathscr{S}(\rr^n)$ for any positive time $t>0$, i.e.,
\begin{equation}\label{p1}
\forall t>0, \quad u(t)=e^{-t\mathscr{L}}u_0 \in \mathscr{S}(\rr^n),
\end{equation}
where $(e^{-t\mathscr{L}})_{t \geq 0}$ denotes the contraction semigroup generated by $\mathscr{L}$. It follows from (\ref{pav4.5}) that the solution to the evolution equation (\ref{p2}) is given by
\begin{equation}\label{bn1}
v(t)=e^{tP}v_0=(\sqrt{\rho})^{-1}e^{-t(\mathscr{L}+\frac{1}{2}\textrm{Tr}(B))}(\sqrt{\rho}v_0)=e^{-\frac{t}{2}\textrm{Tr}(B)}(\sqrt{\rho})^{-1}e^{-t\mathscr{L}}(\sqrt{\rho}v_0),
\end{equation}
for all $t \geq 0$.
By using from (\ref{pav2}) that
$$\sqrt{\rho}\partial_{x_i}\big((\sqrt{\rho})^{-1}u\big)=e^{-\frac{1}{4}\langle Q_{\infty}^{-1}x,x\rangle}\partial_{x_i}(e^{\frac{1}{4}\langle Q_{\infty}^{-1}x,x\rangle}u)=\Big(\partial_{x_i}+\frac{1}{2}(Q_{\infty}^{-1}x)_i\Big)u,$$
where $(Q_{\infty}^{-1}x)_i$ denotes the $i^{\textrm{th}}$ coordinate,
since the matrix $Q_{\infty}^{-1}$ is symmetric, we notice that for all $\alpha=(\alpha_1,...,\alpha_n) \in \mathbb{N}^n$, $\beta=(\beta_1,...,\beta_n) \in \mathbb{N}^n$, $t \geq 0$,
$$x^{\alpha}\partial_{x}^{\beta}v(t)=e^{-\frac{t}{2}\textrm{Tr}(B)}(\sqrt{\rho})^{-1}
x^{\alpha}\Big(\partial_{x_1}+\frac{1}{2}(Q_{\infty}^{-1}x)_1\Big)^{\beta_1}
...\Big(\partial_{x_n}+\frac{1}{2}(Q_{\infty}^{-1}x)_n\Big)^{\beta_n}e^{-t\mathscr{L}}(\sqrt{\rho}v_0).$$
We deduce from (\ref{p1}) that for all $t>0$,
$$x^{\alpha}\Big(\partial_{x_1}+\frac{1}{2}(Q_{\infty}^{-1}x)_1\Big)^{\beta_1}
...\Big(\partial_{x_n}+\frac{1}{2}(Q_{\infty}^{-1}x)_n\Big)^{\beta_n}e^{-t\mathscr{L}}(\sqrt{\rho}v_0) \in L^2,$$
since $\sqrt{\rho}v_0 \in L^2$, because $v_0 \in L_{\mu}^2$. This implies that 
$$\forall t >0,\forall \alpha, \beta \in \mathbb{N}^n, \quad x^{\alpha}\partial_{x}^{\beta}v(t) \in L_{\mu}^2.$$
This ends the proof of Proposition~\ref{prop1}. 
\end{proof}

\subsection{Spectrum of hypoelliptic Ornstein-Ulhenbeck operators}
Let 
$$P=\frac{1}{2}\textrm{Tr}(Q\nabla_x^2)+\langle Bx,\nabla_x\rangle,\quad x \in \rr^n,$$ 
be a hypoelliptic Ornstein-Uhlenbeck operator, which admits the invariant measure $d\mu(x)=\rho(x) dx$. 
The Markov semigroup $(T(t))_{t \geq 0}$ defined in (\ref{gen}) extends to a strongly continuous semigroup of positive contractions in $L^r_{\mu}=L^r(\rr^n,d\mu)$, for every $1 \leq r <+\infty$. 
We denote $(P_r,D_r)$ the generator of $(T(t))_{t \geq 0}$ in the $L_{\mu}^r$ space.

The result by Metafune, Pallara and Priola~\cite{Metafune_al2002} (Theorem~3.1) shows that for all $1<r<+\infty$, the spectrum of the generator $P_r$ is discrete and is composed of eigenvalues with finite algebraic multiplicities given by
\begin{equation}\label{rt1}
\sigma(P_r)=\Big\{\sum_{\lambda \in \sigma(B)}\lambda k_{\lambda} : k_{\lambda} \in \mathbb{N}\Big\},
\end{equation}
where $\mathbb{N}$ stands for the set of non-negative integers.
When $r=1$, the spectrum of the generator $P_1$ in the $L^1_{\mu}$ space is the closed left half-plane 
$$\sigma(P_1)=\{z \in \CC : \textrm{Re }z \leq 0\}.$$
Furthermore, each complex number $z$ with a negative real part $\textrm{Re} z<0$ is an eigenvalue~\cite{Metafune_al2002} (Theorem~5.1).

The following proposition shows that we can recover from~\cite{HPS} (Theorem~1.2.2)
the description of the spectrum of hypoelliptic Ornstein-Ulhenbeck operators on $L^2_{\mu}$ spaces with respect to invariant measures:

\medskip

\begin{proposition}\label{rt2}
Let 
$$P=\frac{1}{2}\emph{\textrm{Tr}}(Q\nabla_x^2)+\langle Bx,\nabla_x\rangle, \quad x \in \rr^n,$$ 
be a hypoelliptic Ornstein-Uhlenbeck operator, which admits the invariant measure $d\mu(x)=\rho(x) dx$. Then, the spectrum of the operator 
$P :  L^2_{\mu} \rightarrow  L^2_{\mu}$ equipped with the domain
\begin{equation}\label{xc5}
D(P)=\{u \in L_{\mu}^2 : Pu \in L_{\mu}^2\},
\end{equation}
is only composed of eigenvalues with finite algebraic multiplicities given by
$$\sigma(P)=\Big\{\sum_{\lambda \in \sigma(B)}\lambda k_{\lambda} : k_{\lambda} \in \mathbb{N}\Big\}.$$
\end{proposition}

\medskip

\begin{proof}
Starting again from the fact that the quadratic operator
$$\mathscr{L}=q^w(x,D_x)=\frac{1}{2}|Q^{\frac{1}{2}}D_x|^2+\frac{1}{8}|Q^{\frac{1}{2}}Q_{\infty}^{-1}x|^2-i\Big\langle\Big(\frac{1}{2}QQ_{\infty}^{-1}+B\Big)x,D_x\Big\rangle,$$
has a Weyl symbol with a zero singular space $S=\{0\}$ and a non-negative real part $\textrm{Re }q \geq 0$,
we deduce from~\cite{HPS} (Theorem~1.2.2) that the spectrum of the operator $\mathscr{L}$ acting on $L^2=L^2(\rr^n,dx)$ equipped with the domain (\ref{dom1}) is only composed of eigenvalues with finite algebraic multiplicities given by
\begin{equation}\label{jk1}
\sigma(\mathscr{L})=\Big\{\sum_{\substack{\lambda \in \sigma(F)\\  -i\lambda \in \CC_+}} (r_{\lambda}+2k_{\lambda})(-i\lambda) : k_{\lambda} \in \mathbb{N}\Big\},
\end{equation}
with $\CC_+=\{z \in \CC : \textrm{Re }z>0\}$,
where $F$ denotes the Hamilton map (\ref{F}) of the Weyl symbol $q$, and where
$r_{\lambda}$ stands for the dimension of the space of generalized eigenvectors of $F$ in $\CC^{2n}$ associated to the eigenvalue $\lambda$. 
By using that the mappings 
$$
\begin{array}{cc}
T :  L^2_{\mu} & \rightarrow  L^2\\
\ \ \ u &  \mapsto  \sqrt{\rho}u
\end{array}, \qquad \begin{array}{cc}
T^{-1} :  L^2 & \rightarrow L^2_{\mu}\\
\quad \ \ \ u & \mapsto \sqrt{\rho}^{-1}u
\end{array},
$$
are isometric, it follows from (\ref{pav4.5}) and (\ref{jk1}) that the spectrum of the operator $P$ acting on $L_{\mu}^2=L^2(\rr^n,d\mu)$ equipped with the domain (\ref{xc5}) is only composed of eigenvalues with finite algebraic multiplicities exactly given by
\begin{equation}\label{jk2}
\sigma(P)=\Big\{-\frac{1}{2}\textrm{Tr}(B)+\sum_{\substack{\lambda \in \sigma(F)\\  -i\lambda \in \CC_+}} (r_{\lambda}+2k_{\lambda})(i\lambda) : k_{\lambda} \in \mathbb{N}\Big\}.
\end{equation}
We notice from (\ref{pav1}) and Corollary~\ref{coro1} that 
\begin{equation}\label{jk3}
\{\lambda \in \sigma(F) :  -i\lambda \in \CC_+\}=\Big\{\pm \frac{i}{2} \mu : \mu \in \sigma(B),\ \pm \mu \in \CC_+\Big\}
=\Big\{- \frac{i}{2} \mu : \mu \in \sigma(B)\Big\},
\end{equation}
since $\sigma(B) \subset \CC_-$. Furthermore, we also deduce from (\ref{pav1}) and Corollary~\ref{coro1} that $r_{\lambda}$ the dimension of the space of generalized eigenvectors of $F$ in $\CC^{2n}$ associated to the eigenvalue $\lambda=- \frac{i}{2} \mu$, with $\mu \in \sigma(B)$, is exactly equal to $\tilde{r}_{\mu}$ the dimension of the space of generalized eigenvectors of $B$ in $\CC^{n}$ associated to the eigenvalue $\mu$, since
$$\sigma\Big(-\frac{i}{2}B\Big) \cap \sigma\Big(\frac{i}{2}B\Big)=\emptyset,$$
because $\sigma(B) \subset \CC_-$.
It follows from (\ref{jk2}) and (\ref{jk3}) that 
$$\sigma(P)=\Big\{-\frac{1}{2}\textrm{Tr}(B)+\frac{1}{2}\sum_{\mu \in \sigma(B)} (\tilde{r}_{\mu}+2k_{\mu})\mu : k_{\mu} \in \mathbb{N}\Big\}=\Big\{\sum_{\mu \in \sigma(B)}\mu k_{\mu} : k_{\mu} \in \mathbb{N}\Big\},$$
since
$$\sum_{\mu \in \sigma(B)}\mu \tilde{r}_{\mu}=\textrm{Tr}(B).$$
This ends the proof of Proposition~\ref{rt2}.
\end{proof}

As noticed in~\cite{Metafune_al2002}, it is interesting to underline that the spectrum of the hypoelliptic Ornstein-Ulhenbeck operator 
$$P=\frac{1}{2}\textrm{Tr}(Q\nabla_x^2)+\langle Bx,\nabla_x\rangle, \quad x \in \rr^n,$$ 
only depends on the spectrum of $B$ the matrix  appearing in the transport part, and not on the diffusion matrix $Q$. We remark, however, that the function space $L^2_{\mu}$ on which acts the operator $P$ depends on $Q$.

\subsection{Global hypoelliptic estimates}
In the work~\cite{lanco}, Bramanti, Cupini, Lanconelli and Priola established global $L^p(\rr^n,dx)$ estimates for the elliptic directions of hypoelliptic Ornstein-Uhlenbeck operators for every $1<p<+\infty$.

More specifically, let  
$$P=\frac{1}{2}\textrm{Tr}(Q\nabla_x^2)+\langle Bx,\nabla_x\rangle, \quad x \in \rr^n,$$ 
be a hypoelliptic Ornstein-Uhlenbeck operator, which does not necessarily admit an invariant measure. 
We define $(V_k)_{k \geq 0}$ the vector subspaces
\begin{equation}\label{jen2.1}
V_k=\big(\textrm{Ran}(Q^{\frac{1}{2}})+\textrm{Ran}(BQ^{\frac{1}{2}})+...+\textrm{Ran}(B^kQ^{\frac{1}{2}})\big) \cap \rr^n \subset \rr^n, \quad k \geq 0,
\end{equation} 
where the notation $\textrm{Ran}$ denotes the range.
Since the Ornstein-Uhlenbeck operator $P$ is hypoelliptic, the Kalman rank condition holds
$$\textrm{Rank}[Q^{\frac{1}{2}},BQ^{\frac{1}{2}},\dots, B^{n-1}Q^{\frac{1}{2}}]=n.$$
We can therefore consider the smallest integer $0 \leq k_0 \leq n-1$ satisfying
\begin{equation}\label{kal2}
\textrm{Rank}[Q^{\frac{1}{2}},BQ^{\frac{1}{2}},\dots, B^{k_0}Q^{\frac{1}{2}}]=n,
\end{equation}
and we notice from (\ref{jen2.1}) and (\ref{kal2}) that 
\begin{equation}\label{byebye}
V_0 \subsetneq V_1 \subsetneq ... \subsetneq V_{k_0}=\rr^n.
\end{equation}
The work by Lanconelli and Polidoro~\cite{lanco2} shows that a fan orthonormal basis
\begin{equation}\label{ui0}
\mathcal{B}=(e_1,...,e_n), \qquad \qquad V_k=\textrm{Span}\{e_j : 1 \leq j \leq \textrm{dim }V_k\}, \quad 0 \leq k \leq k_0,
\end{equation} 
for the subspaces $V_0,...,V_{k_0}$, may be chosen so that the Ornstein-Uhlenbeck operator in these new coordinates can be written as 
\begin{equation}\label{jen1}
\tilde{P}=\sum_{i,j=1}^{p_0}\tilde{a}_{i,j}\partial_{x_i,x_j}^2+\sum_{i,j=1}^nc_{i,j}x_j\partial_{x_i},
\end{equation}
with $1 \leq p_0 \leq n$, where $\tilde{A}=(\tilde{a}_{i,j})_{1 \leq i,j \leq p_0} \in M_{p_0}(\rr)$ is symmetric positive definite and $C=(c_{i,j})_{1 \leq i,j \leq n}$ has the block structure
$$C=\left( \begin{array}{ccccc}
* & * & \cdots & * & * \\
C_1 & * & \cdots & * & * \\
0 & C_2 & \ddots & * & * \\
\vdots & \ddots & \ddots & * & * \\
0 & \cdots & 0 & C_{k_0} & * \\
\end{array} \right)\in M_{n}(\rr),$$
where $C_j$ is a $p_j \times p_{j-1}$ block with rank $p_j$ for all $j=1,...,k_0$ satisfying
$$p_0 \geq p_1 \geq p_2 \geq ... \geq p_{k_0} \geq 1,$$
where 
$$p_0=\textrm{dim }V_0, \quad p_i=\textrm{dim }V_i-\textrm{dim }V_{i-1}, \quad 1 \leq i \leq k_0.$$
Bramanti, Cupini, Lanconelli and Priola established in~\cite{lanco} the following global $L^p$ estimates  
\begin{equation}\label{ellip}
\forall 1<p<+\infty, \exists C_p>0, \forall u \in C_0^{\infty}(\rr^n), \quad \sum_{i,j=1}^{p_0}\|\partial_{x_i,x_j}^2u\|_{L^p} \leq C_p(\|\tilde{P}u\|_{L^p}+\|u\|_{L^p}),
\end{equation}
where $\|\cdot\|_{L^p}$ denotes the $L^p(\rr^n,dx)$ norm.

On the other hand, Farkas and Lunardi provided in~\cite{lunardi1} sharp embeddings for the domains of hypoelliptic Ornstein-Uhlenbeck operators with invariant measures acting on $L^2_{\mu}$ spaces, in some anisotropic Sobolev spaces. In order to recall the results of~\cite{lunardi1}, we consider the sets of indices
\begin{equation}\label{ui2}
I_0=\{j \in \mathbb{N} : 1\leq j \leq \textrm{dim }V_0\}, \quad I_k=\{j \in \mathbb{N} : \textrm{dim }V_{k-1}+1\leq j \leq \textrm{dim }V_k\},
\end{equation}
for $1 \leq k \leq k_0$, providing the partition
$$I_0 \sqcup I_1 \sqcup .... \sqcup I_{k_0}=\{1,...,n\}.$$
The 1-dimensional Hermite polynomials are defined by
$$h_n(x)=\frac{(-1)^n}{\sqrt{n!}}e^{\frac{x^2}{2}}\frac{d^n}{dx^n}e^{-\frac{x^2}{2}}, \quad n \geq 0,\ x \in \rr.$$
By using the fact that the symmetric matrix $Q_{\infty}$ is positive definite, we can introduce an orthogonal matrix $U$, fixed once for all, such that $UQ_{\infty}U^{-1}=\textrm{diag}[\lambda_1,...,\lambda_n]$ is diagonal. We define for any multi-index $\beta=(\beta_1,...,\beta_n) \in \mathbb{N}^n$,
$$H_{\beta}(x)=\prod_{j=1}^nh_{\beta_j}\Big(\frac{(U\mathcal{T}^{-1}x)_{j}}{\sqrt{\lambda_j}}\Big), \quad x \in \rr^n,$$
where $\mathcal{T}$ denotes the invertible matrix representing the change of basis from the canonical basis of $\rr^n$ to the basis $\mathcal{B}$ defined in (\ref{ui0}). 
As eigenfunctions of the selfadjoint non-positive Ornstein-Uhlenbeck operator 
\begin{equation}\label{hj10}
AH_{\beta}=\Big(\frac{1}{2}\textrm{Tr}(Q_{\infty}\nabla_x^2)-\frac{1}{2}\langle x,\nabla_x\rangle\Big)H_{\beta}=-\frac{|\beta|}{2}H_{\beta},
\end{equation}
with $|\beta|=\beta_1+...+\beta_n$, these polynomials constitute an orthonormal basis of $L^2_{\mu}$. For any $s>0$,
the Sobolev space $H^s(\rr^n,d\mu)$ is defined as the domain of the operator $(\sqrt{I-A})^s$, i.e., the set of $L^2_{\mu}$ functions satisfying
\begin{equation}\label{hj12}
\|u\|_{H^s(\rr^n,d\mu)}^2=\|(\sqrt{I-A})^su\|_{L^2_{\mu}}^2
=\sum_{\beta \in \mathbb{N}^n}\Big(1+\frac{|\beta|}{2}\Big)^{s}|(u,H_{\beta})_{L^2_{\mu}}|^2<+\infty.
\end{equation}
For $s_0,s_1,...,s_{k_0}>0$, the anisotropic Sobolev space $H^{s_0,s_1,....,s_{k_0}}(\rr^n,d\mu)$ is defined as the space of $L^2_{\mu}$ functions satisfying
$$\|u\|_{H^{s_0,s_1,....,s_{k_0}}(\rr^n,d\mu)}^2=\sum_{\beta \in \mathbb{N}^n}\sum_{k=0}^{k_0}\Big(1+\sum_{j \in I_k}\frac{\beta_j}{2}\Big)^{s_k}|(u,H_{\beta})_{L_{\mu}^2}|^2,$$ 
where the sets of indices $(I_k)_{0 \leq k \leq k_0}$ are defined in (\ref{ui2}).
The result of~\cite{lunardi1} (Theorem~8) shows that the domain of the infinitesimal generator of the Ornstein-Ulhenbeck semigroup $(T(t))_{t \geq 0}$ satisfies the following embedding into the anisotropic Sobolev space
\begin{equation}\label{ui5}
D(P) \subset H^{2,\frac{2}{3},\frac{2}{5}....,\frac{2}{2k_0+1}}(\rr^n,d\mu) \subset H^{\frac{2}{2k_0+1}}(\rr^n,d\mu).
\end{equation}

The following proposition shows that we can recover the weaker embedding of the domain into the isotropic Sobolev space
$$D(P) \subset H^{\frac{2}{2k_0+1}}(\rr^n,d\mu),$$
from the global $L^2$ subelliptic estimates established for general accretive quadratic operators with zero singular spaces in~\cite{karel} (Theorem~1.2.1). 
This indicates that the general result of~\cite{karel} manages to capture the exact loss of derivatives with respect to the elliptic case 
$$\delta=2-\frac{2}{2k_0+1}=\frac{4k_0}{2k_0+1}>0,$$ 
even if it misses to provide stronger estimates in the less degenerate frequency directions.

\medskip

\begin{proposition}\label{subest}
Let 
$$P=\frac{1}{2}\emph{\textrm{Tr}}(Q\nabla_x^2)+\langle Bx,\nabla_x\rangle, \quad x \in \rr^n,$$ 
be a hypoelliptic Ornstein-Uhlenbeck operator, which admits the invariant measure $d\mu(x)=\rho(x) dx$.
Then, there exists a positive constant $C>0$ such that 
$$\forall v \in D(P), \quad \|v\|_{H^{\frac{2}{2k_0+1}}(\rr^n,d\mu)}  \leq C(\|Pv\|_{L_{\mu}^2}+\|v\|_{L_{\mu}^2}),$$
where $k_0$ denotes the smallest integer $0 \leq k_0 \leq n-1$ satisfying
$$\emph{\textrm{Rank}}[Q^{\frac{1}{2}},BQ^{\frac{1}{2}},\dots, B^{k_0}Q^{\frac{1}{2}}]=n,$$
with $\|\cdot\|_{L_{\mu}^2}$ the norm in the $L^2(\rr^n,d\mu)$ space, where the Sobolev space $H^{\frac{2}{2k_0+1}}(\rr^n,d\mu)$ is defined in (\ref{hj12}).
\end{proposition}

\medskip

\begin{proof}
By using that
$$e^{-\frac{1}{4}\langle Q_{\infty}^{-1}x,x\rangle}\partial_{x_i}(e^{\frac{1}{4}\langle Q_{\infty}^{-1}x,x\rangle}u)=\Big(\partial_{x_i}+\frac{1}{2}(Q_{\infty}^{-1}x)_i\Big)u,$$
where $(Q_{\infty}^{-1}x)_i$ denotes the $i^{\textrm{th}}$ coordinate,
since the matrix $Q_{\infty}^{-1}$ is symmetric, we obtain from (\ref{pav2}) and (\ref{hj10}) that 
\begin{multline}\label{pav11n}
-\sqrt{\rho}A\big((\sqrt{\rho})^{-1}u\big) =-\frac{1}{2}e^{-\frac{1}{4}\langle Q_{\infty}^{-1}x,x\rangle}\Big(\sum_{i,j=1}^{n}\tilde{q}_{i,j}\partial_{x_i,x_j}^2-\sum_{j=1}^nx_j\partial_{x_j}\Big)(e^{\frac{1}{4}\langle Q_{\infty}^{-1}x,x\rangle}u)\\
=-\frac{1}{2}\sum_{i,j=1}^{n}\tilde{q}_{i,j}\Big(\partial_{x_i}+\frac{1}{2}(Q_{\infty}^{-1}x)_i\Big)\Big(\partial_{x_j}+\frac{1}{2}(Q_{\infty}^{-1}x)_j\Big)u+\frac{1}{2}\sum_{j=1}^nx_j\Big(\partial_{x_j}+\frac{1}{2}(Q_{\infty}^{-1}x)_j\Big)u,
\end{multline}
with $Q_{\infty}=(\tilde{q}_{i,j})_{1 \leq i,j \leq n}$.
By using that 
$$\sum_{i,j=1}^{n}\tilde{q}_{i,j}\partial_{x_i}\big((Q_{\infty}^{-1}x)_j\big)=\textrm{Tr}(Q_{\infty}Q_{\infty}^{-1})=n,$$
it follows that 
$$\sqrt{\rho}(I-A)\big((\sqrt{\rho})^{-1}u\big)
=\frac{1}{2}|Q_{\infty}^{\frac{1}{2}}D_x|^2u+\frac{1}{8}|Q_{\infty}^{-\frac{1}{2}}x|^2u+\Big(1-\frac{n}{4}\Big)u.$$
We deduce from (\ref{hj10}) that the family $(H_{\beta}\sqrt{\rho})_{\beta \in \mathbb{N}^n}$ is a Hilbert basis of the $L^2(\rr^n)$ space composed by the eigenvalues of the harmonic oscillator
\begin{equation}\label{hj11}
\mathscr{H}(H_{\beta}\sqrt{\rho})=\Big(\frac{1}{2}|Q_{\infty}^{\frac{1}{2}}D_x|^2+\frac{1}{8}|Q_{\infty}^{-\frac{1}{2}}x|^2\Big)(H_{\beta}\sqrt{\rho})=\Big(\frac{|\beta|}{2}+\frac{n}{4}\Big)H_{\beta}\sqrt{\rho}.
\end{equation}
It follows from (\ref{hj12}) that for all $s>0$,
\begin{multline}\label{hj121}
\|v\|_{H^s(\rr^n,d\mu)}^2=\|(\sqrt{I-A})^sv\|_{L^2_{\mu}}^2
=\sum_{\beta \in \mathbb{N}^n}\Big(1+\frac{|\beta|}{2}\Big)^{s}|(v,H_{\beta})_{L^2_{\mu}}|^2\\
=\sum_{\beta \in \mathbb{N}^n}\Big(1+\frac{|\beta|}{2}\Big)^{s}|(\sqrt{\rho}v,H_{\beta}\sqrt{\rho})_{L^2}|^2=\Big\|\Big(1+\mathscr{H}-\frac{n}{4}\Big)^{\frac{s}{2}}(\sqrt{\rho}v)\Big\|_{L^2}^2,
\end{multline}
where the fractional power of the harmonic oscillator is defined by functional calculus.
It follows from (\ref{hj121}) that for any $s>0$,
\begin{equation}\label{hj122}
\|v\|_{H^s(\rr^n,d\mu)} \sim \|\langle\sqrt{\mathscr{H}}\rangle^s(\sqrt{\rho}v)\|_{L^2}.
\end{equation}
Starting anew from the fact that the quadratic operator
$$\mathscr{L}=q^w(x,D_x)=\frac{1}{2}|Q^{\frac{1}{2}}D_x|^2+\frac{1}{8}|Q^{\frac{1}{2}}Q_{\infty}^{-1}x|^2-i\Big\langle\Big(\frac{1}{2}QQ_{\infty}^{-1}+B\Big)x,D_x\Big\rangle,$$
has a Weyl symbol with a zero singular space $S=\{0\}$ and a non-negative real part $\textrm{Re }q \geq 0$, we can apply the result of~\cite{karel} (Theorem~1.2.1) to show that there exists a positive constant $C>0$ such that 
\begin{equation}\label{xc2}
\forall u \in D(\mathscr{L}), \quad \big\|\textrm{Op}^w\big(\langle (x,\xi) \rangle^{\frac{2}{2k_0+1}}\big)u\big\|_{L^2} \leq C(\|\mathscr{L}u\|_{L^2}+\|u\|_{L^2}),
\end{equation}
where the domain $D(\mathscr{L})$ is defined in (\ref{dom1}) and where $k_0$ denotes the smallest integer $0 \leq k_0 \leq 2n-1$ satisfying
$$\Big(\bigcap_{j=0}^{k_0}\textrm{Ker}
\big[\textrm{Re }F(\textrm{Im }F)^j \big]\Big)\cap \rr^{2n}=\{0\},$$
with $F$ the Hamilton map of the quadratic operator $\mathscr{L}$. We notice from (\ref{pav30}) and (\ref{pav31}) that this integer $k_0$ corresponds exactly to the smallest integer $0 \leq k_0 \leq n-1$ satisfying 
$$\textrm{Rank}[Q^{\frac{1}{2}},BQ^{\frac{1}{2}},\dots, B^{k_0}Q^{\frac{1}{2}}]=n.$$
We easily check the equivalence of the norms
\begin{equation}\label{gh1}
\big\|\textrm{Op}^w\big(\langle (x,\xi) \rangle^{\frac{2}{2k_0+1}}\big)u\big\|_{L^2} \sim \big\|\textrm{Op}^w\big(\big\langle \big((2\sqrt{2})^{-1}Q_{\infty}^{-\frac{1}{2}}x,(\sqrt{2})^{-1}Q_{\infty}^{\frac{1}{2}}\xi\big)\big\rangle^{\frac{2}{2k_0+1}}\big)u\big\|_{L^2}.
\end{equation}
As noticed in~\cite{herau} (Proposition A.4), see also the discussion in~\cite{OPPS} (pp.~4021-4022), the following equivalences of norms hold 
\begin{equation}\label{gh2}
\big\|\textrm{Op}^w\big(\big\langle \big((2\sqrt{2})^{-1}Q_{\infty}^{-\frac{1}{2}}x,(\sqrt{2})^{-1}Q_{\infty}^{\frac{1}{2}}\xi\big)\big\rangle^{\frac{2}{2k_0+1}}\big)u\big\|_{L^2} \sim \|\langle \sqrt{\mathscr{H}} \rangle^{\frac{2}{2k_0+1}}u\|_{L^2}
\end{equation}
and
\begin{equation}\label{gh2f}
\big\|\textrm{Op}^w\big(\langle (x,\xi)\rangle^{\frac{2}{2k_0+1}}\big)u\big\|_{L^2} \sim \|\langle (x,D_x) \rangle^{\frac{2}{2k_0+1}}u\|_{L^2},
\end{equation}
where the harmonic oscillator $\mathscr{H}$ is defined in (\ref{hj11}) and where the operator $\langle (x,D_x) \rangle^{\frac{2}{2k_0+1}}$ is defined by the functional calculus of the operator 
$$\langle (x,D_x) \rangle^2=1+|D_x|^2+|x|^2.$$
By taking $v=(\sqrt{\rho})^{-1}u \in L_{\mu}^2$, with $u \in L^2$, we deduce from (\ref{pav4.5}), (\ref{xc5}), (\ref{hj122}), (\ref{xc2}), (\ref{gh1}) and (\ref{gh2}) that there exists a positive constant $C_0>0$ such that for all $v \in D(P)$,  
$$\|v\|_{H^{\frac{2}{2k_0+1}}(\rr^n,d\mu)} \leq C_0(\|\mathscr{L}u\|_{L^2}+\|u\|_{L^2})
\leq C_0(\|Pv\|_{L_{\mu}^2}+C_0\Big(1+\frac{1}{2}|\textrm{Tr}(B)|\Big)\|v\|_{L_{\mu}^2}.$$
This ends the proof of Proposition~\ref{subest}.
\end{proof}

\subsection{Spectral instabilities and resolvent estimates}
Let 
\begin{equation}\label{yu1}
P=\frac{1}{2}\textrm{Tr}(Q\nabla_x^2)+\langle Bx,\nabla_x\rangle, \quad x \in \rr^n,
\end{equation} 
be a hypoelliptic Ornstein-Uhlenbeck operator, which admits the invariant measure $d\mu(x)=\rho(x) dx$. 
Recalling that $B \in M_n(\rr)$ and thus $\textrm{Tr}(B) \in \rr$, we notice from (\ref{pav4.5}) and (\ref{pav5}) that the operator $P$ is selfadjoint on $L^2_{\mu}$ if and only if the operator $\mathscr{L}$ is selfadjoint on $L^2$, that is, when  
$$\frac{1}{2}Q_{\infty}^{-1}Q+B^T=0.$$
When this condition does not hold, the operator is non-selfadjoint. In this case, it is well-known that the resolvent of such an operator may take very large values in norm far away from its spectrum and that the spectrum may be very unstable under small perturbations. These two features are highly related. Indeed, studying the level lines of the norm of the resolvent of an operator provides substantial information about its spectral stability. We recall from~\cite{roch} that when $A$ is a closed unbounded linear operator with a dense domain on a complex Hilbert space $H$, its $\varepsilon$-pseudospectrum
$$\sigma_{\varepsilon}(A)=\Big\{z \in \cc:  \ \|(A-z)^{-1}\| \geq \frac{1}{\varepsilon} \Big\}, \quad \eps>0,$$
with the convention that $\|(A-z)^{-1}\|=+\infty$ for any point $z \in \sigma(A)$, can be defined in an equivalent way in term of the spectrum of perturbations of the operator
$$\sigma_{\eps}(A)=\bigcup_{B \in \mathcal{L}(H), \ \|B\|_{\mathcal{L}(H)} \leq \eps}{\sigma(A+B)},$$
where $\mathcal{L}(H)$ stands for the set of bounded linear operators on $H$. There exists an extensive literature on the notion of pseudospectrum. We refer the reader to~\cite{trefethen,trefethen2} and the references therein for a detailed account of this topic.

The study of pseudospectrum is non-trivial only for non-selfadjoint operators, or more precisely for non-normal operators. Indeed, the classical formula 
\begin{equation}\label{1}
\forall z \not\in \sigma(A), \quad \|(A-z)^{-1}\| = \frac{1}{\textrm{dist}(z,\sigma(A))}, 
\end{equation}
implies that the resolvent of a normal operator 
\begin{equation}\label{normal}
\forall z \not\in \sigma(A), \forall \zeta \not\in \sigma(A^*), \quad (A-z)^{-1}(A^*-\zeta)^{-1}=(A^*-\zeta)^{-1}(A-z)^{-1},
\end{equation}
cannot blow up far away from its spectrum~\cite{kato} (Chap.~V, Sect.~3.8, formula (3.31)). It ensures the stability of the spectrum under small perturbations 
\begin{equation}\label{2}
\sigma_{\eps}(A) = \{z \in \cc :\  \textrm{dist}(z,\sigma(A)) \leq \eps\}. 
\end{equation}
However, the formula (\ref{1}) does not hold anymore for non-normal operators and the resolvent of such operators  may become very large in norm far away from the spectrum. It implies that the spectra of these operators may become very unstable under small perturbations.

Consider now again the hypoelliptic Ornstein-Uhlenbeck operator (\ref{yu1}). We observe from (\ref{pav4.5}) that the operator $P$ commutes with its $L^2_{\mu}$-adjoint $P^*$, 
\begin{equation}\label{commu1}
\forall u \in \mathscr{S}(\rr^n), \quad [P,P^*](\sqrt{\rho}^{-1}u)=0,
\end{equation} 
when acting on $\sqrt{\rho}^{-1}\mathscr{S}(\rr^n)$,
if and only if the operator $\mathscr{L}$ commutes with its $L^2$-adjoint~$\mathscr{L}^*$,  
\begin{equation}\label{commu2}
\forall u \in \mathscr{S}(\rr^n), \quad [\mathscr{L},\mathscr{L}^*]u=0,
\end{equation} 
when acting on $\mathscr{S}(\rr^n)$, since $[\mathscr{L},\mathscr{L}^*]=-\sqrt{\rho}[P,P^*]\sqrt{\rho}^{-1}$. We deduce from Lemma~\ref{normal8} that the condition (\ref{commu2}) holds
if and only if the following commutator is zero 
$$[Q_{\infty}^{-1}Q,B^T]=0.$$ 
Furthermore, Lemma~\ref{normal8} indicates that, when the condition (\ref{commu1}) does not hold, then the operator $P$ is non-normal, and that the operator $P$ can be normal only if it is elliptic, i.e. when the symmetric matrix $Q$ is positive definite. When the operator $P$ is elliptic, both situations can occur. The operator $P$ is for instance selfadjoint, thus normal, when $B=-I_n$, whereas it is non-normal when
$$Q=\left(\begin{array}{cc}
1 & 0\\
0 & 1
\end{array}\right), \quad B=\left(\begin{array}{cc}
-1 & 1\\
0 & -1
\end{array}\right).$$

The study of resolvent estimates for hypoelliptic Ornstein-Uhlenbeck operators is therefore non-trivial. The resolvents of these operators may exhibit rapid growth in norm far away from their spectra. The following theorem points out that this is actually the case for any non-normal elliptic Ornstein-Uhlenbeck operator satisfying
$$[Q_{\infty}^{-1}Q,B^T] \neq 0.$$
More specifically, we deduce from the results in~\cite{duke} that the resolvent of such an operator exhibits rapid growth in norm along all the half-lines contained in the particular open angular sector   
$$\Lambda=-\frac{1}{2}\textrm{Tr}(B)-\stackrel{\circ}{\Sigma(q)} \subset \cc,$$
starting from $-\frac{1}{2}\textrm{Tr}(B)$,
where $\stackrel{\circ}{\Sigma(q)}$ denotes the interior of the numerical range $\Sigma(q)=\overline{q(\rr^{2n})}$ of the Weyl symbol (\ref{pav6}).

\medskip

\begin{theorem}\label{esti}
Let 
$$P=\frac{1}{2}\emph{\textrm{Tr}}(Q\nabla_x^2)+\langle Bx,\nabla_x\rangle, \quad x \in \rr^n,$$ 
be an elliptic Ornstein-Uhlenbeck operator, i.e., when $Q$ is a symmetric positive definite matrix, which admits the invariant measure $d\mu(x)=\rho(x) dx$ and satisfies the condition $[Q_{\infty}^{-1}Q,B^T]\neq 0$, implying that $P$ is non-normal. Setting
$$m_-=\inf_{(x,\xi) \in \mathcal{B}}\Big\langle\Big(-\frac{1}{2}QQ_{\infty}^{-1}-B\Big)x,\xi\Big\rangle, \quad m_+=\sup_{(x,\xi) \in \mathcal{B}}\Big\langle\Big(-\frac{1}{2}QQ_{\infty}^{-1}-B\Big)x,\xi\Big\rangle,$$
with
$$\mathcal{B}=\Big\{(x,\xi) \in \rr^{2n} :\ \frac{1}{2}|Q^{\frac{1}{2}}\xi|^2+\frac{1}{8}|Q^{\frac{1}{2}}Q_{\infty}^{-1}x|^2=1\Big\},$$
then the resolvent blows-up in norm along all the half-lines 
$$\forall m_- < \lambda < m_+, \forall N \in \mathbb{N}, \quad \lim_{t \to +\infty}\frac{\big\|\big(P+\frac{1}{2}\emph{\textrm{Tr}}(B)+t(1+i\lambda)\big)^{-1}\big\|_{\mathcal{L}(L^2_{\mu})}}{t^N}=+\infty,$$
where $\|\cdot\|_{\mathcal{L}(L_{\mu}^2)}$ denotes the operator norm in the space of bounded operators on $L^2_{\mu}=L^2(\rr^n,d\mu)$.
\end{theorem}

\medskip

This result highlights that the resolvent of a non-normal elliptic Ornstein-Uhlenbeck operator as in Theorem~\ref{esti},
blows up rapidly in norm far away from its spectrum 
$$\sigma(P)=\Big\{\sum_{\lambda \in \sigma(B)}\lambda k_{\lambda} : k_{\lambda} \in \mathbb{N}\Big\},$$
described in Proposition~\ref{rt2}.

\begin{proof}
Under the assumptions of the theorem, we observe from Lemma~\ref{normal8} that the quadratic operator on $L^2$,
$$\mathscr{L}=q^w(x,D_x)=\frac{1}{2}|Q^{\frac{1}{2}}D_x|^2+\frac{1}{8}|Q^{\frac{1}{2}}Q_{\infty}^{-1}x|^2-i\Big\langle\Big(\frac{1}{2}QQ_{\infty}^{-1}+B\Big)x,D_x\Big\rangle,$$
is elliptic and does not commute with its $L^2$-adjoint $\mathscr{L}^*$,
$$[\mathscr{L},\mathscr{L}^*]\neq 0,$$
when acting on $\mathscr{S}(\rr^n)$. 
The interior of the numerical range $\Sigma(q)=\overline{q(\rr^{2n})}$ of its Weyl symbol (\ref{pav6}) is given by
$$ \stackrel{\circ}{\Sigma(q)}=\{z=t(1+i\lambda) \in \cc^* : \  m_- < \lambda < m_+, \ t>0\},$$
where the constants $m_{\pm}$ are defined in the statement of Theorem~\ref{esti}. It follows from~\cite{duke} (Theorem~2.2.1) that for all $m_- < \lambda < m_+$ and $N \in \mathbb{N}$, there exist $h_0>0$ and a family of the Schwartz functions $(u_h)_{0<h \leq h_0}$ satisfying
\begin{equation}\label{gh0}
\forall 0<h \leq h_0, \quad \|u_h\|_{L^2}=1, \qquad \|q^w(x,hD_x)u_h-(1+i\lambda)u_h\|_{L^2}=O(h^N),
\end{equation}
when $h \to 0$. We notice that
\begin{equation}\label{gh1r}
\forall t>0, \quad \tilde{T}_t\big(q^w(x,D_x)-t(1+i\lambda)\big)\tilde{T}_t^{-1}=t\big(q^w(x,t^{-1}D_x)-(1+i\lambda)\big),
\end{equation}
where $\tilde{T}_t$ stands for the isometry of $L^2(\rr^n)$ defined by $(\tilde{T}_tu)(x)=t^{\frac{n}{4}}u(\sqrt{t}x)$. We deduce from (\ref{gh0}) and (\ref{gh1r}) that for all $m_- < \lambda < m_+$, $N \in \mathbb{N}$, 
\begin{multline}\label{df1}
\lim_{t \to +\infty}\frac{\big\|\big(\mathscr{L}-t(1+i\lambda)\big)^{-1}\big\|_{\mathcal{L}(L^2)}}{t^N}\\ =\lim_{t \to +\infty}\frac{\big\|\big(q^w(x,t^{-1}D_x)-(1+i\lambda)\big)^{-1}\big\|_{\mathcal{L}(L^2)}}{t^{N+1}}=+\infty,
\end{multline}
where $\|\cdot\|_{\mathcal{L}(L^2)}$ denotes the operator norm in the space of bounded operators on $L^2$.
By using that the mappings 
$$
\begin{array}{cc}
T :  L^2_{\mu} & \rightarrow  L^2\\
\ \ \ u &  \mapsto  \sqrt{\rho}u
\end{array}, \qquad \begin{array}{cc}
T^{-1} :  L^2 & \rightarrow L^2_{\mu}\\
\quad \ \ \ u & \mapsto \sqrt{\rho}^{-1}u
\end{array},
$$
are isometric, we notice from (\ref{pav4.5}) that 
\begin{equation}\label{ul1}
\big\|\big(\mathscr{L}-t(1+i\lambda)\big)^{-1}\big\|_{\mathcal{L}(L^2)}=\Big\|\Big(P+\frac{1}{2}\textrm{Tr}(B)+t(1+i\lambda)\Big)^{-1}\Big\|_{\mathcal{L}(L_{\mu}^2)},
\end{equation}
where $\|\cdot\|_{\mathcal{L}(L_{\mu}^2)}$ denotes the operator norm in the space of bounded operators on $L^2_{\mu}=L^2(\rr^n,d\mu)$.
The statement of Theorem~\ref{esti} then directly follows from (\ref{df1}) and (\ref{ul1}).
\end{proof}

The result of~\cite{duke} does not apply directly to the degenerate case when the symmetric matrix~$Q$ is not positive definite, but the same blow-up phenomena actually occur for all degenerate hypoelliptic Ornstein-Uhlenbeck operators, i.e. when the diffusion matrix $Q$ is degenerate. The following result shows that the resolvent of any degenerate hypoelliptic Ornstein-Uhlenbeck operator exhibits rapid growth in norm along all the half-lines contained in the open half-plane
$$-\frac{1}{2}\textrm{Tr}(B)+\cc_- \subset \cc,$$
starting from $-\frac{1}{2}\textrm{Tr}(B)$.

\medskip

\begin{theorem}\label{esti9}
Let 
$$P=\frac{1}{2}\emph{\textrm{Tr}}(Q\nabla_x^2)+\langle Bx,\nabla_x\rangle, \quad x \in \rr^n,$$ 
be a degenerate hypoelliptic Ornstein-Uhlenbeck operator, i.e. when $\det(Q)=0$, which admits the invariant measure $d\mu(x)=\rho(x) dx$.
Then, the resolvent blows-up in norm along all the half-lines 
$$\forall z \in \cc_+, \forall N \in \mathbb{N}, \quad \lim_{t \to +\infty}\frac{\big\|\big(P+\frac{1}{2}\emph{\textrm{Tr}}(B)+tz\big)^{-1}\big\|_{\mathcal{L}(L^2_{\mu})}}{t^N}=+\infty,$$
where $\|\cdot\|_{\mathcal{L}(L_{\mu}^2)}$ denotes the operator norm in the space of bounded operators on $L^2_{\mu}=L^2(\rr^n,d\mu)$.
\end{theorem}

\medskip

\begin{proof}
Since $Q$ is degenerate, we deduce from Lemma~\ref{jk1b} that the numerical range $\Sigma(q)=\overline{q(\rr^{2n})}$ of the quadratic symbol (\ref{pav6})
is equal to 
$$\Sigma(q)=\{z \in \cc :  \textrm{Re }z \geq 0\}.$$
Furthermore, for all $z \in \cc_+$, there exists $(x_0,\xi_0) \in \rr^{2n}$ such that 
\begin{equation}\label{xiao1}
z=q(x_0,\xi_0), \qquad \{\textrm{Re }q,\textrm{Im }q\}(x_0,\xi_0) <0.
\end{equation}
Then, the result of Zworski~\cite{zworski1} (p. 2956) and~\cite{zworski2} (Theorem, p. 300) allows us to construct semiclassical quasimodes for the operator $q^w(x,hD_x)$ at any point $z \in \cc_+$. More specifically, for all $z \in \cc_+$ and $N \in \mathbb{N}$, we can find $h_0>0$ and a family of the Schwartz functions $(u_h)_{0<h \leq h_0}$ satisfying
\begin{equation}\label{xiao2}
\forall 0<h \leq h_0, \quad \|u_h\|_{L^2}=1, \qquad \|q^w(x,hD_x)u_h-zu_h\|_{L^2}=O(h^N),
\end{equation}
when $h \to 0$. The proof of Theorem~\ref{esti9} is then completed by following the very same lines as in the proof of Theorem~\ref{esti}.
\end{proof}

\begin{figure}[ht]
\caption{The resolvent $\|(P-z)^{-1}\|_{\mathcal{L}(L^2_{\mu})}$ blows-up along all the half-lines contained in the dark region starting from $-\frac{1}{2}\textrm{Tr}(B)$. The spectrum is represented with black dots. The left figure corresponds to the non-normal elliptic case as in Theorem~\ref{esti} and the right one to the degenerate case.}
\begin{minipage}[b]{0.3\linewidth}
\centering
\includegraphics[scale=0.58]{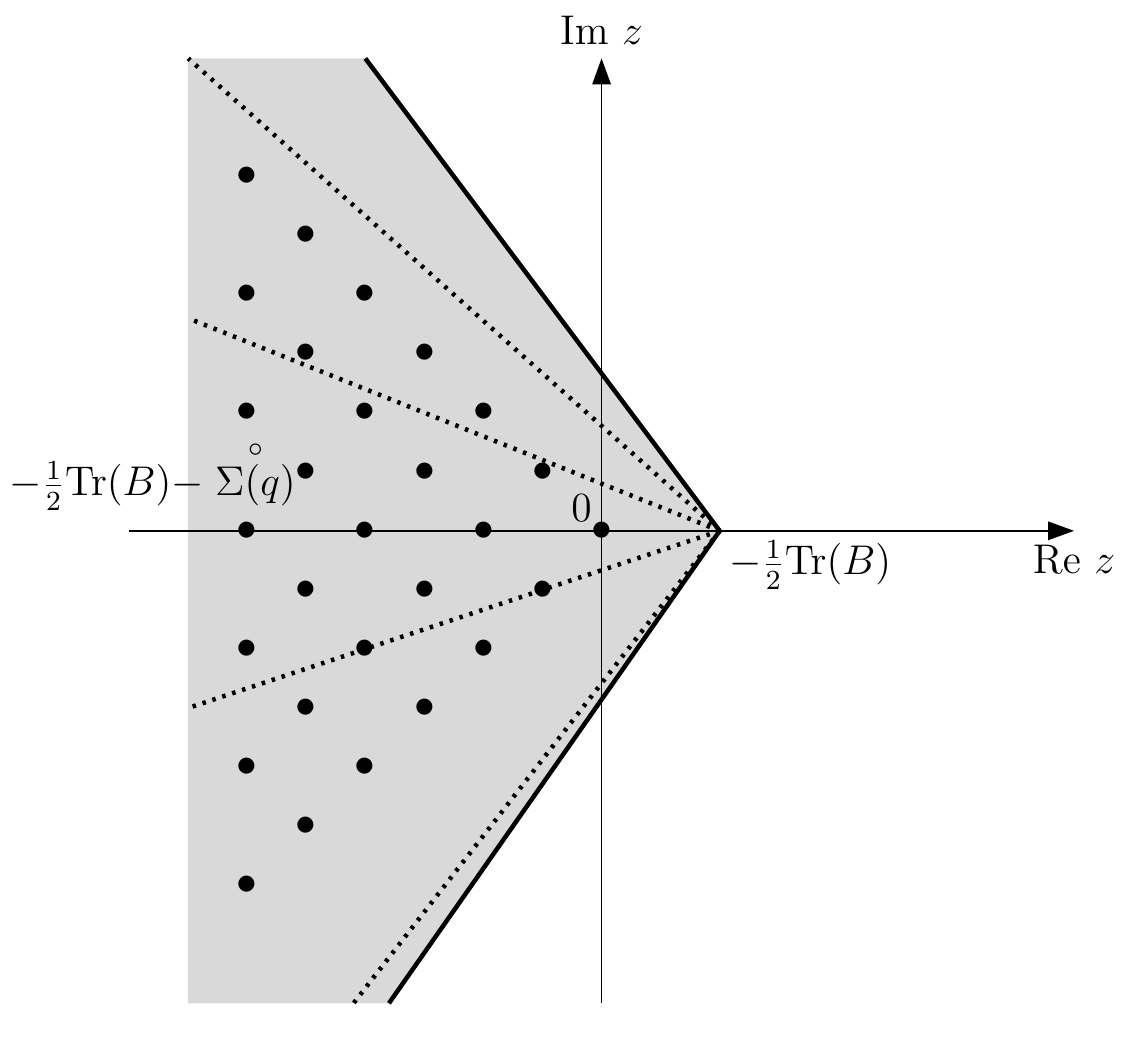}
\label{fig:figure1.1}
\end{minipage}
\hspace{2cm}
\begin{minipage}[b]{0.3\linewidth}
\centering
\includegraphics[scale=0.58]{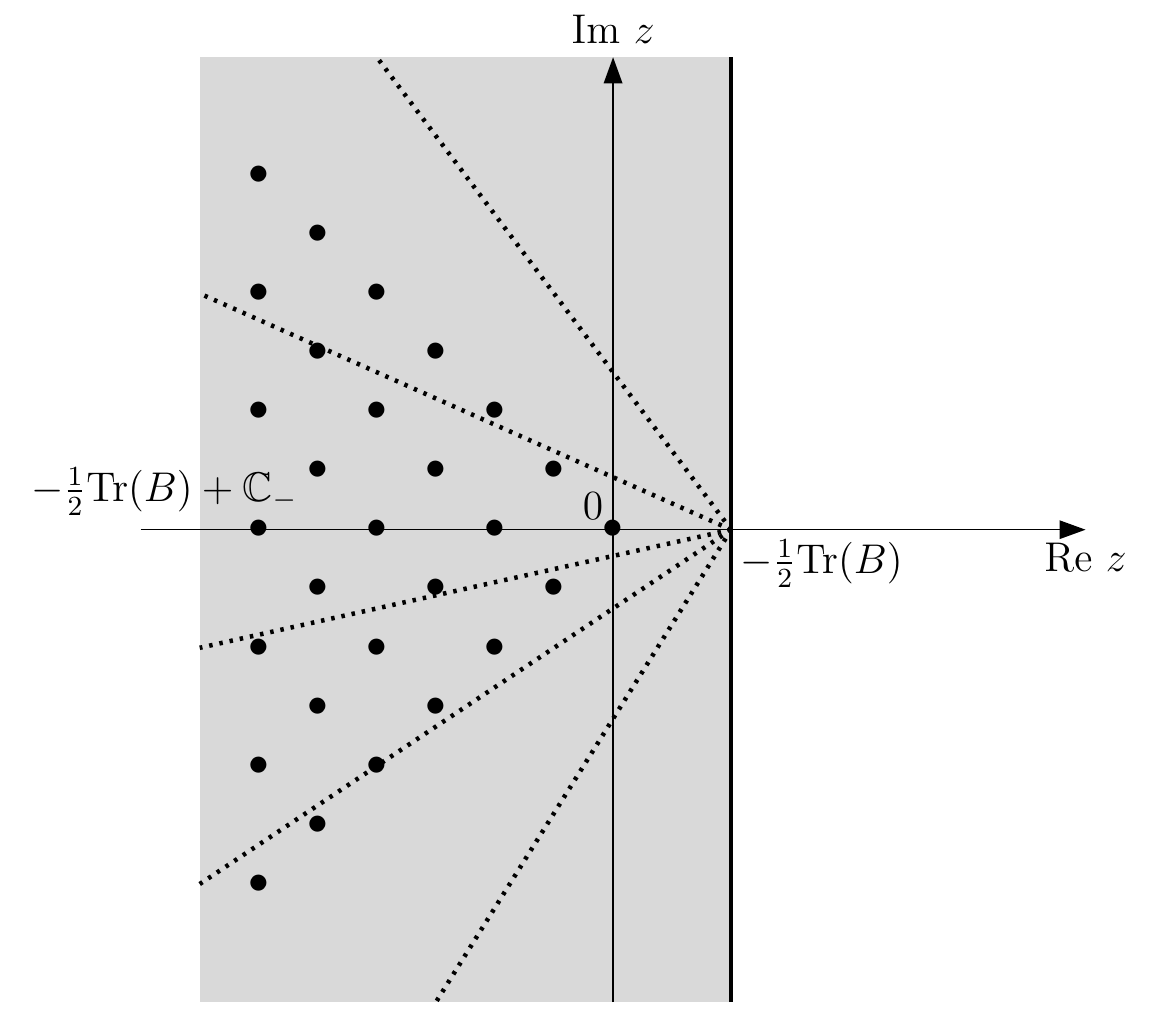}
\label{fig:figure2.1}
\end{minipage}
\hspace{0.5cm}
\end{figure}

\medskip

\begin{remark} The above result for the construction of semiclassical quasimodes is a particular case of a general existence result of semiclassical quasimodes for pseudodifferential operators violating the condition $(\overline{\Psi})$, i.e., for operators whose adjoints violate the Nirenberg-Treves condition, also called condition $(\Psi)$, see Definition 26.4.6 and Theorem~26.4.12 in \cite{hormander3}. This result indicates that, when the principal symbol of $P_0-z$ violates the condition $(\overline{\Psi})$, then there exist some semiclassical quasimodes of the type (\ref{xiao2}) for the operator $P^w(x,hD_x;h)-z$.
The condition (\ref{xiao1}) corresponds to the particular case when the principal symbol violates the condition $(\overline{\Psi})$ by changing sign in the prohibited sense at the first order, whereas the general result holds true more generally without any restriction on the order, finite or infinite, at which the prohibited change of sign given by the violation of the condition $(\overline{\Psi})$ is done. The existence of this result was first mentioned in~\cite{dencker}.  It is an adaptation to the semiclassical setting of the quasimode construction intiated by Moyer~\cite{moyer} and completed in all dimensions by H\"ormander~\cite{hormander3} (Theorem 26.4.7).  A complete proof of this adaptation in the semiclassical setting is given in \cite{these}. 
\end{remark}

\medskip

Despite the blow-up phenomena described in Theorems~\ref{esti} and~\ref{esti9}, it was shown in~\cite{OPPS} that the subelliptic properties enjoyed 
by quadratic operators with zero singular spaces imply that their resolvents must stay bounded in certain unbounded regions of the resolvent set with a specific geometry. The following theorem makes explicit how these results apply for hypoelliptic Ornstein-Uhlenbeck operators with invariant measures:  

\bigskip

\begin{theorem}\label{esti1}
Let 
$$P=\frac{1}{2}\emph{\textrm{Tr}}(Q\nabla_x^2)+\langle Bx,\nabla_x\rangle, \quad x \in \rr^n,$$ 
be a hypoelliptic Ornstein-Uhlenbeck operator, which admits the invariant measure $d\mu(x)=\rho(x) dx$. Then, there exist some positive constants $C,c>0$ such that the resolvent estimate
$$\forall z \in \Gamma_{k_0}, \quad  \|(P-z)^{-1}\|_{\mathcal{L}(L_{\mu}^2)} \leq C \Big|z-\Big(1-\frac{1}{2}\emph{\textrm{Tr}}(B)\Big)\Big|^{-\frac{1}{2k_0+1}},$$
holds in the subset of the resolvent set $\Gamma_{k_0}  \subset \cc \setminus \sigma(P)$ defined as
$$\Gamma_{k_0}=\Big\{z \in \CC : \ \emph{\textrm{Re }}z \leq  \frac{1}{2}\big(1-\emph{\textrm{Tr}}(B)\big), \  \Big|\emph{\textrm{Re }}z -\Big(1-\frac{1}{2}\emph{\textrm{Tr}}(B)\Big)\Big|  \leq c\Big|z-\Big(1-\frac{1}{2}\emph{\textrm{Tr}}(B)\Big)\Big|^{\frac{1}{2k_0+1}}\Big\},$$
with $\|\cdot\|_{\mathcal{L}(L_{\mu}^2)}$ the operator norm in the space of bounded operators on $L^2_{\mu}=L^2(\rr^n,d\mu)$, 
where $k_0$ denotes the smallest integer $0 \leq k_0 \leq n-1$ satisfying
$$\emph{\textrm{Rank}}[Q^{\frac{1}{2}},BQ^{\frac{1}{2}},\dots, B^{k_0}Q^{\frac{1}{2}}]=n.$$
Moreover, we have
$$\forall z \in \cc,\ \emph{\textrm{Re }}z>-\frac{1}{2}\emph{\textrm{Tr}}(B)>0, \quad  \|(P-z)^{-1}\|_{\mathcal{L}(L_{\mu}^2)} \leq \frac{1}{\emph{\textrm{Re }}z+\frac{1}{2}\emph{\textrm{Tr}}(B)}.$$
\end{theorem}

\bigskip

\begin{proof}
Following the proof of Proposition~\ref{subest}, we deduce from (\ref{xc2}) and (\ref{gh2f}) that the quadratic operator (\ref{pav5}) enjoys the following global subelliptic estimate
$$\exists C>0, \forall u \in D(\mathscr{L}), \quad \|\langle(x,D_x)\rangle^{\frac{2}{2k_0+1}}u\|_{L^2} \leq C(\|\mathscr{L}u\|_{L^2}+\|u\|_{L^2}),$$
where $k_0$ denotes the smallest integer $0 \leq k_0 \leq n-1$ satisfying
$$\textrm{Rank}[Q^{\frac{1}{2}},BQ^{\frac{1}{2}},\dots, B^{k_0}Q^{\frac{1}{2}}]=n.$$
As explained in~\cite{OPPS} (pp. 4021-4023), this estimate is shown to extend as 
\begin{equation}\label{est3}
\exists C>0, \forall u \in D(\mathscr{L}), \forall \nu \in \rr, \quad  \|\langle(x,D_x)\rangle^{\frac{2}{2k_0+1}} u\|_{L^2} \leq C(\|\mathscr{L}u-i\nu u\|_{L^2}+\|u\|_{L^2}).
\end{equation}
Using this estimate and standard functional analytic estimates we derive the localization of the spectrum
$$\exists c,C>0, \quad \Big\{z \in \cc :\ \textrm{Re }z \geq -\frac{1}{2}, \ \textrm{Re }z +1 \leq c|z+1|^{\frac{1}{2k_0+1}}\Big\} \cap \sigma(\mathscr{L})=\emptyset,$$
together with the resolvent estimate
\begin{equation}\label{dl3}
\|(\mathscr{L}-z)^{-1}\|_{\mathcal{L}(L^2)} \leq C |z+1|^{-\frac{1}{2k_0+1}},
\end{equation}
for all $z \in \cc$ satisfying $\textrm{Re }z \geq -\frac{1}{2}$, $\textrm{Re }z +1 \leq c|z+1|^{\frac{1}{2k_0+1}}$. On the other hand, it follows from (\ref{pav5})
that for all $\textrm{Re }z < 0$,
$$\textrm{Re}(\mathscr{L}u-zu,u)_{L^2}=\frac{1}{2}\|Q^{\frac{1}{2}}D_xu\|_{L^2}^2+\frac{1}{8}\|Q^{\frac{1}{2}}Q_{\infty}^{-1}xu\|_{L^2}^2-\textrm{Re }z\|u\|_{L^2}^2 \geq |\textrm{Re }z|\|u\|_{L^2}^2.$$
We deduce from (\ref{jk1}) and the Cauchy-Schwarz inequality that for all $\textrm{Re }z < 0$,
\begin{equation}\label{fgh}
\|(\mathscr{L}-z)^{-1}\|_{\mathcal{L}(L^2)} \leq \frac{1}{|\textrm{Re }z|}.
\end{equation}
It follows from (\ref{pav4.5}), (\ref{ul1}) and (\ref{dl3}) that the resolvent estimate
$$\forall z \in \tilde{\Gamma}_{k_0}, \quad  \Big\|\Big(P+\frac{1}{2}\textrm{Tr}(B)-z\Big)^{-1}\Big\|_{\mathcal{L}(L_{\mu}^2)} \leq C |z-1|^{-\frac{1}{2k_0+1}},$$
holds in the following subset of the complex plane
$$\tilde{\Gamma}_{k_0}=\Big\{z \in \CC : \ \textrm{Re }z \leq  \frac{1}{2}, \  |\textrm{Re }z -1|  \leq c|z-1|^{\frac{1}{2k_0+1}}\Big\} \subset \cc \setminus \sigma\Big(P+\frac{1}{2}\textrm{Tr}(B)\Big).$$
Moreover, we observe from (\ref{fgh}) that 
$$\forall z \in \cc,\ \textrm{Re }z>0, \quad  \Big\|\Big(P+\frac{1}{2}\textrm{Tr}(B)-z\Big)^{-1}\Big\|_{\mathcal{L}(L_{\mu}^2)} \leq \frac{1}{\textrm{Re }z}.$$
This ends the proof of Theorem~\ref{esti1}.
\end{proof}

\begin{figure}[ht]
\caption{Spectrum of the hypoelliptic Ornstein-Uhlenbeck operator $P$ acting on $L_{\mu}^2$ and subset $\Gamma_{k_0}$.}
\centerline{\includegraphics[scale=0.8]{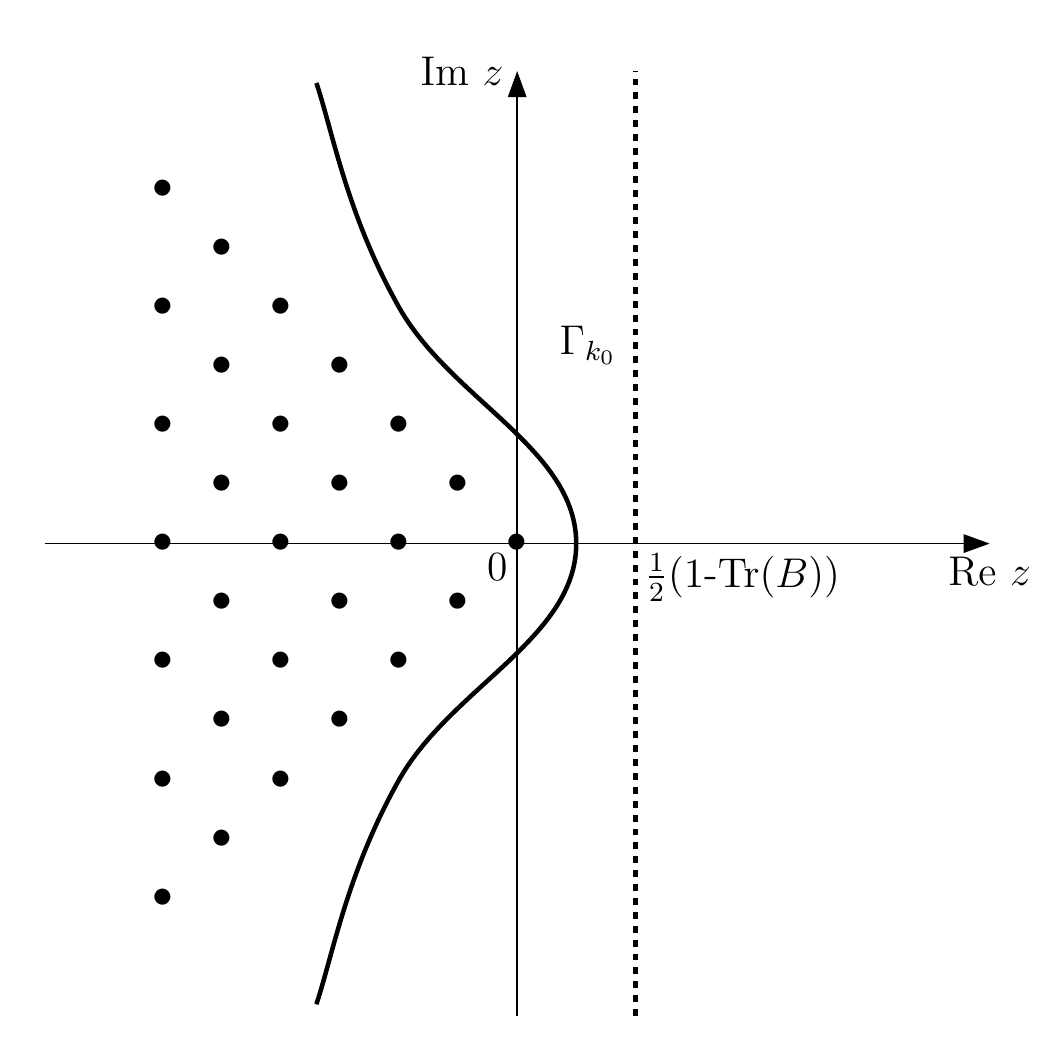}}
\end{figure}  

As explained in~\cite{OPPS} (pp. 4023-4025), the resolvent estimates as in Theorem~\ref{esti1} are key for establishing results on exponential return to equilibrium for the Ornstein-Uhlenbeck semigroup.

\subsection{Exponential return to equilibrium for hypoelliptic Ornstein-Uhlenbeck operators}

The following result establishes the result of exponential return to equilibrium for hypoelliptic Ornstein-Uhlenbeck operators with invariant measures:

\begin{theorem}\label{return}
Let 
$$P=\frac{1}{2}\emph{\textrm{Tr}}(Q\nabla_x^2)+\langle Bx,\nabla_x\rangle, \quad x \in \rr^n,$$ 
be a hypoelliptic Ornstein-Uhlenbeck operator, which admits the invariant measure $d\mu(x)=\rho(x) dx$.
Setting
$$\tau_0=\inf_{\mu \in \sigma(B)}\emph{\textrm{Re}}(-\mu)=-\sup_{\mu \in \sigma(B)}\emph{\textrm{Re }}\mu>0,$$
then for all $0 \leq \tau<\tau_0$, there exists a positive constant $C>0$ such that 
$$\forall t \geq 0, \forall v \in L_{\mu}^2, \quad 
\Big\|e^{tP}v-\Big(\int_{\rr^{n}}v(x)\rho(x)dx\Big)\Big\|_{L_{\mu}^2} \leq Ce^{-\tau t}\|v\|_{L_{\mu}^2}.$$
\end{theorem}

\medskip

\begin{proof}
We begin by noticing from (\ref{pav1}) that the parameter $\tau_0>0$ is positive.
We remark that the quadratic operator (\ref{pav5}) is real
$$\mathscr{L}=-\frac{1}{2}|Q^{\frac{1}{2}}\nabla_x|^2+\frac{1}{8}|Q^{\frac{1}{2}}Q_{\infty}^{-1}x|^2-\Big\langle\Big(\frac{1}{2}QQ_{\infty}^{-1}+B\Big)x,\nabla_x\Big\rangle,$$ 
in the sense that $\mathscr{L}u$ is a real-valued function whenever $u$ is a real-valued function since $B$, $Q$ and $Q_{\infty}^{-1} \in M_n(\rr)$. 
Starting now from the fact that the quadratic operator $\mathscr{L}$ 
has a Weyl symbol with a zero singular space $S=\{0\}$ and a non-negative real part, we notice that these properties hold true as well for its $L^2(\rr^n)$-adjoint (\ref{hj1}), since $\mathscr{L}^*$ is a quadratic operator 
whose Weyl symbol $\overline{q}$ is the complex conjugate of the Weyl symbol of~$\mathscr{L}$, and  
whose Hamilton map $\overline{F}$ is the complex conjugate of the Hamilton map~$F$ of~$\mathscr{L}$. 
We deduce from (\ref{pav1}), (\ref{pav4.5}), Proposition~\ref{rt2}, \cite{OPPS} (Theorem~2.1) and~\cite{kato} (Chap.~III, Sect.~6.6, Thm~6.22) that 
$$-\frac{1}{2}\textrm{Tr}(B),$$ 
is the eigenvalue with the lowest real part for both the operators $\mathscr{L}$ and $\mathscr{L}^*$ on $L^2(\rr^n)$. 
Furthermore, this eigenvalue has algebraic multiplicity 1 for both operators.
It therefore follows from (\ref{io0}) that 
\begin{equation}\label{tl1}
\textrm{Ker}\Big(\mathscr{L}+\frac{1}{2}\textrm{Tr}(B)\Big)=\textrm{Ker}\Big(\mathscr{L}^*+\frac{1}{2}\textrm{Tr}(B)\Big)=\CC \sqrt{\rho}.
\end{equation}
Setting
\begin{equation}\label{rty}
\tau_0=2\inf_{\substack{\lambda \in \sigma(F) \\ \textrm{Im }\lambda>0}}\textrm{Im }\lambda>0,
\end{equation}
we deduce from~\cite{OPPS} (Theorem~2.3) that for all $0 \leq \tau<\tau_0$, there exists a positive constant $C>0$ such that for all $t \geq 0$, $u \in L^2(\rr^n)$, 
\begin{equation}\label{bn2}
 \Big\|e^{-t(\mathscr{L}+\frac{1}{2}\textrm{Tr}(B))}u-\frac{(u,\sqrt{\rho})_{L^2}}{\|\sqrt{\rho}\|_{L^2}^2}\sqrt{\rho}\Big\|_{L^2} \leq Ce^{-\tau t}\|u\|_{L^2}.
\end{equation}
We notice from Corollary~\ref{coro1} and (\ref{pav1}) that 
$$\tau_0=\inf_{\mu \in \sigma(B)}\textrm{Re}(-\mu)=-\sup_{\mu \in \sigma(B)}\textrm{Re }\mu>0.$$
On the other hand, we deduce from (\ref{pav2}) that 
\begin{equation}\label{tl2}
\|\sqrt{\rho}\|_{L^2}^2=\int_{\rr^n}\rho(x)dx=\int_{\rr^n}\frac{1}{(2\pi)^{\frac{n}{2}}\sqrt{\det Q_{\infty}}}e^{-\frac{1}{2}\langle Q_{\infty}^{-1}x,x\rangle}dx=1.
\end{equation}
By taking $u=\sqrt{\rho}v \in L^2$ with $v \in L_{\mu}^2$, it follows from (\ref{bn1}) and (\ref{bn2}) that 
\begin{equation}\label{bn3}
\forall t \geq 0, \forall v \in L_{\mu}^2, \quad \Big\|e^{tP}v-\Big(\int_{\rr^{n}}v(x)\rho(x)dx\Big)\Big\|_{L_{\mu}^2} \leq Ce^{-\tau t}\|v\|_{L_{\mu}^2}.
\end{equation}
This ends the proof of Theorem~\ref{return}.
\end{proof}

\subsection{Relative entropy decay for degenerate Fokker-Planck operators}
We consider the Fokker-Planck operator
\begin{equation}\label{jen077}
\mathscr{P}=\frac{1}{2}\textrm{Tr}(Q\nabla_x^2)-\langle Bx,\nabla_x\rangle-\textrm{Tr}(B), \quad x \in \rr^n,
\end{equation}
where $Q=(q_{i,j})_{1 \leq i,j \leq n}$ and $B=(b_{i,j})_{1 \leq i,j \leq n}$ are real $n \times n$-matrices, with $Q$ symmetric positive semidefinite. We assume that the Kalman rank condition 
and the localization of the spectrum of~$B$,
\begin{equation}\label{kal177}
\textrm{Rank}[Q^{\frac{1}{2}},BQ^{\frac{1}{2}},\dots, B^{n-1}Q^{\frac{1}{2}}]=n, \qquad \sigma(B) \subset \mathbb{C}_-,
\end{equation} 
hold. As before, we consider 
\begin{equation}\label{pav277}
\rho(x)=\frac{1}{(2\pi)^{\frac{n}{2}}\sqrt{\det Q_{\infty}}}e^{-\frac{1}{2}\langle Q_{\infty}^{-1}x,x\rangle},
\end{equation}
with 
\begin{equation}\label{pav377}
Q_{\infty}=\int_0^{+\infty}e^{sB}Qe^{sB^T}ds.
\end{equation}
We aim at studying the operator $\mathscr{P}$ acting on the $L^2_{1/\rho}=L^2(\rr^n,\rho(x)^{-1}dx)$ space.
In the recent work~\cite{anton}, Arnold and Erb studied the degenerate parabolic Fokker-Planck equation
$$\partial_tf=\frac{1}{2}\textrm{Tr}(Q\nabla_x^2)f-\langle Bx,\nabla_x\rangle f-\textrm{Tr}(B)f.$$
By employing a new entropy method based on a modified, non-degenerate entropy dissipation like functional,
they established results on exponential convergence of the solution of the Fokker-Planck equation to equilibrium and the exponential decay in relative entropy (logarithmic till quadratic) with a sharp rate. The following proposition shows that the result of~\cite{HPS} for quadratic operators with zero singular spaces enables us to recover the description of the spectrum established in~\cite[Thm~5.3]{anton}: 

\medskip

\begin{proposition}\label{rt277}
Let 
$$\mathscr{P}=\frac{1}{2}\emph{\textrm{Tr}}(Q\nabla_x^2)-\langle Bx,\nabla_x\rangle-\emph{\textrm{Tr}}(B), \quad x \in \rr^n,$$
be a Fokker-Planck operator satisfying the assumption (\ref{kal177}). Then, the spectrum of the operator 
$\mathscr{P} : L^2_{1/\rho} \rightarrow L^2_{1/\rho}$ equipped with the domain
\begin{equation}\label{xc577}
D(\mathscr{P})=\big\{f \in  L^2_{1/\rho} : \mathscr{P}f \in  L^2_{1/\rho}\big\},
\end{equation}
is only composed of eigenvalues with finite algebraic multiplicities exactly given by
$$\sigma(\mathscr{P})=\Big\{\sum_{\lambda \in \sigma(B)}\lambda k_{\lambda} : k_{\lambda} \in \mathbb{N}\Big\}.$$
\end{proposition}

\medskip

\begin{proof}
We begin by associating to the operator $\mathscr{P}$ acting on $L^2_{1/\rho}$, the quadratic operator $\mathfrak{L}$ acting on $L^2$ given by
\begin{equation}\label{pav4.577}
\mathfrak{L}u=-\sqrt{\rho}^{-1}\mathscr{P}\big(\sqrt{\rho}u\big)-\frac{1}{2}\textrm{Tr}(B)u.
\end{equation}
By using the fact that
$$e^{\frac{1}{4}\langle Q_{\infty}^{-1}x,x\rangle}\partial_{x_i}(e^{-\frac{1}{4}\langle Q_{\infty}^{-1}x,x\rangle}u)=\Big(\partial_{x_i}-\frac{1}{2}(Q_{\infty}^{-1}x)_i\Big)u,$$
where $(Q_{\infty}^{-1}x)_i$ denotes the $i^{\textrm{th}}$ coordinate,
since the matrix $Q_{\infty}^{-1}$ is symmetric, we obtain from (\ref{jen077}), (\ref{pav277}) and (\ref{pav4.577}) that 
\begin{multline}\label{pav1177}
\mathfrak{L}u=e^{\frac{1}{4}\langle Q_{\infty}^{-1}x,x\rangle}\Big(-\frac{1}{2}\sum_{i,j=1}^{n}q_{i,j}\partial_{x_i,x_j}^2+\sum_{i,j=1}^nb_{i,j}x_j\partial_{x_i}\Big)(e^{-\frac{1}{4}\langle Q_{\infty}^{-1}x,x\rangle}u)+\frac{1}{2}\textrm{Tr}(B)u=\\
-\frac{1}{2}\sum_{i,j=1}^{n}q_{i,j}\Big(\partial_{x_i}-\frac{1}{2}(Q_{\infty}^{-1}x)_i\Big)\Big(\partial_{x_j}-\frac{1}{2}(Q_{\infty}^{-1}x)_j\Big)u+\sum_{i,j=1}^nb_{i,j}x_j\Big(\partial_{x_i}-\frac{1}{2}(Q_{\infty}^{-1}x)_i\Big)u+\frac{1}{2}\textrm{Tr}(B)u.
\end{multline}
It follows from (\ref{yp1}) that 
\begin{multline*}
\mathfrak{L}
=-\frac{1}{2}\langle Q\nabla_x,\nabla_x\rangle-\frac{1}{8}\langle QQ_{\infty}^{-1}x,Q_{\infty}^{-1}x\rangle-\frac{1}{2}\langle Bx,Q_{\infty}^{-1}x\rangle\\ +\langle Bx,\nabla_x\rangle+\frac{1}{2}\langle QQ_{\infty}^{-1}x,\nabla_x\rangle+\frac{1}{2}\textrm{Tr}(B)+\frac{1}{4}\textrm{Tr}(QQ_{\infty}^{-1}).
\end{multline*}
With $D_x=i^{-1}\nabla_x$, we deduce from (\ref{not}), (\ref{pav13}) and (\ref{pav14}) that
\begin{equation}\label{pav1277}
\mathfrak{L}=\frac{1}{2}|Q^{\frac{1}{2}}D_x|^2+\frac{1}{8}|Q^{\frac{1}{2}}Q_{\infty}^{-1}x|^2+i\Big\langle \Big(\frac{1}{2} QQ_{\infty}^{-1}+B\Big)x,D_x\Big\rangle.
\end{equation}
We notice that this operator is equal to the $L^2(\rr^n)$-adjoint operator of $\mathscr{L}$
defined in (\ref{hj1}). It follows from (\ref{pav4.5}), Proposition~\ref{rt2} and~\cite{kato} (Chap.~III, Sect.~6.6, Thm~6.22) that 
the spectrum of the operator 
$\mathfrak{L} : L^2 \rightarrow L^2$ equipped with the domain
$$D(\mathfrak{L})=\{u \in L^2 : \mathfrak{L}u \in L^2\},$$
is only composed of eigenvalues with finite algebraic multiplicities exactly given by
\begin{equation}\label{jk377}
\sigma(\mathfrak{L})=\Big\{-\sum_{\lambda \in \sigma(B)}\overline{\lambda} k_{\lambda}-\frac{1}{2}\textrm{Tr}(B) : k_{\lambda} \in \mathbb{N}\Big\},
\end{equation}
since $B \in M_n(\rr)$.
By using the fact that 
$$\mu \in \sigma(B) \Longleftrightarrow \overline{\mu} \in \sigma(B),$$
because $B \in M_n(\rr)$,
it follows from (\ref{pav4.577}) and (\ref{jk377}) that 
$$\sigma(\mathscr{P})=\Big\{\sum_{\lambda \in \sigma(B)}\lambda k_{\lambda} : k_{\lambda} \in \mathbb{N}\Big\}.$$
This ends the proof of Proposition~\ref{rt277}.
\end{proof}

The following proposition and its corollary show that the result of~\cite{OPPS} for quadratic operators with zero singular spaces enables us to recover the exponential decay in quadratic relative entropy established in~\cite[Thm.~4.6]{anton}:

\medskip

\begin{proposition}\label{return77}
Let 
$$\mathscr{P}=\frac{1}{2}\emph{\textrm{Tr}}(Q\nabla_x^2)-\langle Bx,\nabla_x\rangle-\emph{\textrm{Tr}}(B), \quad x \in \rr^n,$$
be a Fokker-Planck operator satisfying the assumption~(\ref{kal177}). 
Setting
$$\tau_0=\inf_{\mu \in \sigma(B)}\emph{\textrm{Re}}(-\mu)=-\sup_{\mu \in \sigma(B)}\emph{\textrm{Re }}\mu>0,$$
then for all $0 \leq \tau<\tau_0$, there exists a positive constant $C>0$ such that 
$$\forall t \geq 0, \forall f \in L^2_{1/\rho}, \quad 
\Big\|e^{t\mathscr{P}}f-\Big(\int_{\rr^{n}}f(x)dx\Big)\rho\Big\|_{L^2_{1/\rho}} \leq Ce^{-\tau t}\|f\|_{L^2_{1/\rho}},$$
where $(e^{t\mathscr{P}}f)_{t \geq 0}$ denotes the semigroup in $L^2_{1/\rho}$ associated to $\mathscr{P}$.
\end{proposition}

\medskip

\begin{proof}
We begin by noticing from (\ref{kal177}) that the parameter $\tau_0>0$ is positive.
We remark that the quadratic operator (\ref{pav1277}) is real
$$\mathfrak{L}=-\frac{1}{2}|Q^{\frac{1}{2}}\nabla_x|^2+\frac{1}{8}|Q^{\frac{1}{2}}Q_{\infty}^{-1}x|^2+\Big\langle \Big(\frac{1}{2} QQ_{\infty}^{-1}+B\Big)x,\nabla_x\Big\rangle,$$ 
since $B$, $Q$ and $Q_{\infty}^{-1} \in M_n(\rr)$. 
As noticed in the proofs of Theorem~\ref{return} and Proposition~\ref{rt277}, the quadratic operator $\mathfrak{L}=\mathscr{L}^*$ has a Weyl symbol with a zero singular space and a non-negative real part. These properties hold true as well for its $L^2(\rr^n)$-adjoint $(\mathfrak{L})^*=\mathscr{L}$, where $\mathscr{L}$ is the quadratic operator (\ref{pav5}).
We deduce from (\ref{tl1}) that 
$$-\frac{1}{2}\textrm{Tr}(B),$$ 
is the eigenvalue with the lowest real part for both the operators $\mathfrak{L}$ and $\mathfrak{L}^*$ on $L^2$,
since $\mathfrak{L}=\mathscr{L}^*$ and $\mathfrak{L}^*=\mathscr{L}$. Furthermore, this eigenvalue has algebraic multiplicity 1 for both operators $\mathfrak{L}$ and $\mathfrak{L}^*$.
It follows from (\ref{tl1}) that 
$$\textrm{Ker}\Big(\mathfrak{L}+\frac{1}{2}\textrm{Tr}(B)\Big)=\textrm{Ker}\Big(\mathfrak{L}^*+\frac{1}{2}\textrm{Tr}(B)\Big)=\CC \sqrt{\rho}.$$
The operator $\mathfrak{L}$ is accretive and generates a contraction semigroup $(e^{-t\mathfrak{L}})_{t \geq 0}$ on $L^2$.
Setting
\begin{equation}\label{rty77}
\tau_0=2\inf_{\substack{\lambda \in \sigma(\overline{F}) \\ \textrm{Im }\lambda>0}}\textrm{Im }\lambda>0,
\end{equation}
we deduce from~\cite{OPPS} (Theorem~2.3) that for all $0 \leq \tau<\tau_0$, there exists a positive constant $C>0$ such that for all $t \geq 0$, $u \in L^2$, 
\begin{equation}\label{bn277}
 \Big\|e^{-t(\mathfrak{L}+\frac{1}{2}\textrm{Tr}(B))}u-\frac{(u,\sqrt{\rho})_{L^2}}{\|\sqrt{\rho}\|_{L^2}^2}\sqrt{\rho}\Big\|_{L^2} \leq Ce^{-\tau t}\|u\|_{L^2},
\end{equation}
since $\overline{F}$ is the Hamilton map of the quadratic operator $\mathfrak{L}$.
We notice from Corollary~\ref{coro1} and (\ref{kal177}) that 
$$\tau_0=\inf_{\mu \in \sigma(B)}\textrm{Re}(-\mu)=-\sup_{\mu \in \sigma(B)}\textrm{Re }\mu>0.$$
On the other hand, the solution to the Cauchy problem
\begin{equation}\label{p277}
\left\lbrace\begin{array}{c}
\partial_tv=\mathscr{P}v,\\
v|_{t=0}=f \in L_{1/\rho}^2,
\end{array}\right.
\end{equation}
is given by
\begin{equation}\label{bn177}
\forall t \geq 0, \quad v(t)=e^{t\mathscr{P}}f=\sqrt{\rho}e^{-t(\mathfrak{L}+\frac{1}{2}\textrm{Tr}(B))}(\sqrt{\rho}^{-1}f),
\end{equation}
where $(e^{-t\mathfrak{L}})_{t \geq 0}$ denotes the contraction semigroup generated by $\mathfrak{L}$.
By taking $u=\sqrt{\rho}^{-1}f \in L^2$ with $f \in L_{1/\rho}^2$, it follows from (\ref{tl2}), (\ref{bn277}) and (\ref{bn177}) that 
\begin{equation}\label{bn377}
\forall t \geq 0, \forall f \in L_{1/\rho}^2, \quad \Big\|e^{t\mathscr{P}}f-\Big(\int_{\rr^{n}}f(x)dx\Big)\rho\Big\|_{L_{1/\rho}^2} \leq Ce^{-\tau t}\|f\|_{L_{1/\rho}^2},
\end{equation}
since the mappings 
$$\begin{array}{cc}
\mathfrak{T} :  L^2  & \rightarrow  L^2_{1/\rho}\\
\ u   & \mapsto  \sqrt{\rho}u
\end{array}, \quad \begin{array}{cc}
\mathfrak{T}^{-1} :  L^2_{1/\rho} & \rightarrow L^2\\
\quad \ u & \mapsto \sqrt{\rho}^{-1}u
\end{array},$$
are isometric.
This ends the proof of Proposition~\ref{return77}.
\end{proof}

We directly deduce from Proposition~\ref{return77} the exponential decay of the quadratic relative entropy
$$\mathfrak{e}_2(f|\rho)=\int_{\rr^n}\Big|\frac{f(x)}{\rho(x)}-1\Big|^2\rho(x)dx.$$

\medskip

\begin{corollary}\label{return776}
Let 
$$\mathscr{P}=\frac{1}{2}\emph{\textrm{Tr}}(Q\nabla_x^2)-\langle Bx,\nabla_x\rangle-\emph{\textrm{Tr}}(B), \quad x \in \rr^n,$$
be a Fokker-Planck operator satisfying the assumption (\ref{kal177}). 
Setting
$$\tau_0=\inf_{\mu \in \sigma(B)}\emph{\textrm{Re}}(-\mu)=-\sup_{\mu \in \sigma(B)}\emph{\textrm{Re }}\mu>0,$$
then for all $0 \leq \tau<\tau_0$, there exists a positive constant $C>0$ such that 
\begin{multline*}
\forall t \geq 0, \forall f \in L^2_{1/\rho}, \  \int_{\rr^n}f(x)dx=1 ,\\ 
\mathfrak{e}_2(e^{t\mathscr{P}}f|\rho)=\int_{\rr^n}\Big|\frac{(e^{t\mathscr{P}}f)(x)}{\rho(x)}-1\Big|^2\rho(x)dx \leq Ce^{-2\tau t}\int_{\rr^n}\frac{|f(x)|^2}{\rho(x)}dx,
\end{multline*}
where $(e^{t\mathscr{P}}f)_{t \geq 0}$ denotes the semigroup in $L^2_{1/\rho}$ associated to $\mathscr{P}$.
\end{corollary}

\section{Elements of proof}\label{ele}

The main purpose of this section is to check that the quadratic operator (\ref{pav4.5}) has a zero singular space. We also compute explicitly the eigenvalues of the Hamilton map (\ref{F}), and provide a necessary and sufficient condition for the quadratic operator to be normal.

\subsection{Computation of the quadratic operator $\mathscr{L}$}\label{pav4}
We begin by noticing from (\ref{pav0}) and (\ref{pav3}) that for all $t \geq 0$,
\begin{equation}\label{pav0.1}
Q_{\infty}=Q_t+e^{tB}Q_{\infty}e^{tB^T}.
\end{equation}
We deduce from (\ref{pav0.1}) the steady state variance equation
\begin{equation}\label{pav10}
\frac{d}{dt}(Q_t+e^{tB}Q_{\infty}e^{tB^T})|_{t=0}=Q+BQ_{\infty}+Q_{\infty}B^T=0.
\end{equation}
We use now
$$e^{-\frac{1}{4}\langle Q_{\infty}^{-1}x,x\rangle}\partial_{x_i}(e^{\frac{1}{4}\langle Q_{\infty}^{-1}x,x\rangle}u)=\Big(\partial_{x_i}+\frac{1}{2}(Q_{\infty}^{-1}x)_i\Big)u,$$
where $(Q_{\infty}^{-1}x)_i$ denotes the $i^{\textrm{th}}$ coordinate
since the matrix $Q_{\infty}^{-1}$ is symmetric, to obtain from (\ref{jen0}), (\ref{pav2}) and (\ref{pav4.5}) that 
\begin{multline*}
\mathscr{L}u = -e^{-\frac{1}{4}\langle Q_{\infty}^{-1}x,x\rangle}\Big(\frac{1}{2}\sum_{i,j=1}^{n}q_{i,j}\partial_{x_i,x_j}^2+\sum_{i,j=1}^nb_{i,j}x_j\partial_{x_i}\Big)(e^{\frac{1}{4}\langle Q_{\infty}^{-1}x,x\rangle}u)-\frac{1}{2}\textrm{Tr}(B)u = \\  
-\frac{1}{2}\sum_{i,j=1}^{n}q_{i,j}\Big(\partial_{x_i}+\frac{1}{2}(Q_{\infty}^{-1}x)_i\Big)\Big(\partial_{x_j}+\frac{1}{2}(Q_{\infty}^{-1}x)_j\Big)u -\sum_{i,j=1}^nb_{i,j}x_j\Big(\partial_{x_i}+\frac{1}{2}(Q_{\infty}^{-1}x)_i\Big)u-\frac{1}{2}\textrm{Tr}(B)u.
\end{multline*}
Furthermore, since
\begin{equation}\label{yp1}
\sum_{i,j=1}^{n}q_{i,j}\partial_{x_i}\big((Q_{\infty}^{-1}x)_j\big)=\textrm{Tr}(QQ_{\infty}^{-1}),
\end{equation}
it follows that 
\begin{multline*}
\mathscr{L} =-\frac{1}{2}\langle Q\nabla_x,\nabla_x\rangle-\frac{1}{8}\langle QQ_{\infty}^{-1}x,Q_{\infty}^{-1}x\rangle-\frac{1}{2}\langle Bx,Q_{\infty}^{-1}x\rangle
\\ -\langle Bx,\nabla_x\rangle-\frac{1}{2}\langle QQ_{\infty}^{-1}x,\nabla_x\rangle-\frac{1}{2}\textrm{Tr}(B)-\frac{1}{4}\textrm{Tr}(QQ_{\infty}^{-1}).
\end{multline*}
By recalling the notation (\ref{not}), we obtain that 
\begin{multline}\label{pav12}
\mathscr{L} = \frac{1}{2}\langle QD_x,D_x\rangle-\frac{1}{8}\langle QQ_{\infty}^{-1}x,Q_{\infty}^{-1}x\rangle-\frac{1}{2}\langle Bx,Q_{\infty}^{-1}x\rangle \\  -i\langle Bx,D_x\rangle-\frac{i}{2}\langle QQ_{\infty}^{-1}x,D_x\rangle-\frac{1}{2}\textrm{Tr}(B)-\frac{1}{4}\textrm{Tr}(QQ_{\infty}^{-1}),
\end{multline}
with $D_x=i^{-1}\nabla_x$.
We deduce from (\ref{pav10}) that
\begin{equation}\label{pav13}
\textrm{Tr}(QQ_{\infty}^{-1})=-\textrm{Tr}(B+Q_{\infty}B^TQ_{\infty}^{-1})=-\textrm{Tr}(B)-\textrm{Tr}(B^T)=-2\textrm{Tr}(B).
\end{equation}
On the other hand, it also follows from (\ref{pav10}) that
\begin{equation}\label{pav14}
-\langle Q_{\infty}^{-1}QQ_{\infty}^{-1}x,x\rangle=\langle Q_{\infty}^{-1}Bx,x\rangle+\langle B^TQ_{\infty}^{-1}x,x\rangle=2\langle Bx,Q_{\infty}^{-1}x\rangle,
\end{equation}
since $Q_{\infty}$ is symmetric. It follows from (\ref{pav12}), (\ref{pav13}) and (\ref{pav14}) that 
\begin{equation}\label{pav15}
\mathscr{L}=\frac{1}{2}|Q^{\frac{1}{2}}D_x|^2+\frac{1}{8}|Q^{\frac{1}{2}}Q_{\infty}^{-1}x|^2-i\Big\langle\Big(\frac{1}{2}QQ_{\infty}^{-1}+B\Big)x,D_x\Big\rangle.
\end{equation}
We notice that the $L^2(\rr^n)$-adjoint operator of $\mathscr{L}$ is given by
\begin{equation}\label{hj1}
\mathscr{L}^*=\frac{1}{2}|Q^{\frac{1}{2}}D_x|^2+\frac{1}{8}|Q^{\frac{1}{2}}Q_{\infty}^{-1}x|^2+i\Big\langle\Big(\frac{1}{2}QQ_{\infty}^{-1}+B\Big)x,D_x\Big\rangle,
\end{equation}
since $B, Q, Q_{\infty} \in M_n(\rr)$, and
$$\textrm{Tr}\Big(\frac{1}{2}QQ_{\infty}^{-1}+B\Big)=0,$$
according to (\ref{pav13}). We check that 
\begin{equation}\label{io0}
\Big(\mathscr{L}+\frac{1}{2}\textrm{Tr}(B)\Big)(\sqrt{\rho})=\Big(\mathscr{L}^*+\frac{1}{2}\textrm{Tr}(B)\Big)(\sqrt{\rho})=0.
\end{equation}
Indeed, it follows from (\ref{pav2}), (\ref{pav10}), (\ref{yp1}) and (\ref{pav13}) that 
\begin{multline}\label{io1}
\Big(\frac{1}{2}|Q^{\frac{1}{2}}D_x|^2\Big)(\sqrt{\rho})=\sum_{1\leq i,j \leq n}q_{i,j}\Big(\frac{1}{4}\partial_{x_j}\big((Q_{\infty}^{-1}x)_i\big)-\frac{1}{8}(Q_{\infty}^{-1}x)_i(Q_{\infty}^{-1}x)_j\Big)\sqrt{\rho}\\
=\Big(\frac{1}{4}\textrm{Tr}(QQ_{\infty}^{-1})-\frac{1}{8}|Q^{\frac{1}{2}}Q_{\infty}^{-1}x|^2\Big)\sqrt{\rho}=\Big(-\frac{1}{2}\textrm{Tr}(B)-\frac{1}{8}|Q^{\frac{1}{2}}Q_{\infty}^{-1}x|^2\Big)\sqrt{\rho}
\end{multline}
and
\begin{equation}\label{io2}
\Big\langle\Big(\frac{1}{2}QQ_{\infty}^{-1}+B\Big)x,D_x\Big\rangle(\sqrt{\rho})=\frac{i}{2}\Big\langle\Big(\frac{1}{2}QQ_{\infty}^{-1}+B\Big)x,Q_{\infty}^{-1}x\Big\rangle \sqrt{\rho}=0,
\end{equation}
since
$$\langle QQ_{\infty}^{-1}x,Q_{\infty}^{-1}x\rangle=\langle Q_{\infty}^{-1}QQ_{\infty}^{-1}x,x\rangle=-\langle Q_{\infty}^{-1}Bx,x\rangle-\langle B^TQ_{\infty}^{-1}x,x\rangle=-2\langle Bx,Q_{\infty}^{-1}x\rangle,$$
because $Q_{\infty}$ is symmetric. The formulas (\ref{io0}) directly follow from (\ref{io1}) and (\ref{io2}).

\subsection{Weyl symbol and Hamilton map of the quadratic operator $\mathscr{L}$}\label{weyl}
By using that the Weyl quantization of the quadratic symbol $x^{\alpha} \xi^{\beta}$, with $(\alpha,\beta) \in \N^{2n}$, $|\alpha+\beta|=2$, is
$$(x^{\alpha} \xi^{\beta})^w=\frac{x^{\alpha}D_x^{\beta}+D_x^{\beta} x^{\alpha}}{2},$$
we observe that the Weyl symbol of the operator
\begin{equation}\label{pav16}
\mathscr{L}=q^w(x,D_x),
\end{equation}
is the quadratic symbol
\begin{equation}\label{pav17}
q(x,\xi)=\frac{1}{2}|Q^{\frac{1}{2}}\xi|^2+\frac{1}{8}|Q^{\frac{1}{2}}Q_{\infty}^{-1}x|^2-i\Big\langle\Big(\frac{1}{2}QQ_{\infty}^{-1}+B\Big)x,\xi\Big\rangle,
\end{equation}
since according to (\ref{pav13}), we have
$$\textrm{Tr}\Big(\frac{1}{2}QQ_{\infty}^{-1}+B\Big)=0.$$
This agrees with the formula~\eqref{pav6}. 
Notice that the polarized form associated to the quadratic form (\ref{pav17}) is
\begin{multline}\label{pav18}
q((x,\xi);(y,\eta)) = \frac{1}{2}\big(q(x+y,\xi+\eta)-q(x,\xi)-q(y,\eta)\big) \\ =
\frac{1}{2}\langle Q^{\frac{1}{2}}\xi,Q^{\frac{1}{2}}\eta \rangle +\frac{1}{8}\langle Q^{\frac{1}{2}}Q_{\infty}^{-1}x,Q^{\frac{1}{2}}Q_{\infty}^{-1}y\rangle 
-\frac{i}{4}\langle(QQ_{\infty}^{-1}+2B)x,\eta\rangle-\frac{i}{4}\langle(QQ_{\infty}^{-1}+2B)y,\xi\rangle. 
\end{multline}
Writing $(\tilde{y},\tilde{\eta})=F(y,\eta)$, we deduce from the definition of the Hamilton map
$$q((x,\xi);(y,\eta))=\sigma((x,\xi),F(y,\eta))=\sigma((x,\xi),(\tilde{y},\tilde{\eta}))=\langle \xi,\tilde{y}\rangle-\langle x,\tilde{\eta}\rangle,$$
that 
\begin{equation}\label{pav19}
(\tilde{y},\tilde{\eta})=F(y,\eta)=\Big(\frac{1}{2}Q\eta-\frac{i}{4}(QQ_{\infty}^{-1}+2B)y,-\frac{1}{8}Q_{\infty}^{-1}QQ_{\infty}^{-1}y
+\frac{i}{4}(QQ_{\infty}^{-1}+2B)^T\eta\Big),
\end{equation}
since $Q$ and $Q_{\infty}$ are symmetric. It follows that 
\begin{equation}\label{pav20}
F=\lv\begin{array}{cc}
 -\frac{i}{4}(QQ_{\infty}^{-1}+2B)&\frac{1}{2}Q \\
 -\frac{1}{8}Q_{\infty}^{-1}QQ_{\infty}^{-1}  & \frac{i}{4}(QQ_{\infty}^{-1}+2B)^T       
\end{array}\rv,
\end{equation}
\begin{equation}\label{pav21}
\textrm{Re }F=\frac{1}{8}\lv\begin{array}{cc}
 0&4Q \\
 -Q_{\infty}^{-1}QQ_{\infty}^{-1}  & 0    
\end{array}\rv
\end{equation}
and
\begin{equation}\label{pav22}
\textrm{Im }F=\frac{1}{4}\lv\begin{array}{cc}
 -(QQ_{\infty}^{-1}+2B)&0 \\
0  & (QQ_{\infty}^{-1}+2B)^T       
\end{array}\rv,
\end{equation}
since $Q$, $Q_{\infty}$ and $B \in M_n(\rr)$.

\subsection{Singular space of the quadratic operator $\mathscr{L}$}\label{singul}

We deduce from (\ref{pav21}) and (\ref{pav22}) that for all $k \geq 0$,
\begin{multline}\label{pav23}
\Big(\bigcap_{j=0}^{k}\textrm{Ker}
\big[\textrm{Re }F(\textrm{Im }F)^j \big]\Big)\cap \rr^{2n}\\ =\big\{(x,\xi) \in \rr^{2n} : \forall 0 \leq j \leq k,\ QQ_{\infty}^{-1}(QQ_{\infty}^{-1}+2B)^jx=0,\ Q(Q_{\infty}^{-1}Q+2B^T)^j\xi=0\big\},
\end{multline}
since $Q$ and $Q_{\infty}$ are symmetric.
We need the following lemma:

\medskip

\begin{lemma}\label{l1}
We have for all $k \geq 0$,
\begin{multline*}
\Big(\bigcap_{j=0}^{k}\emph{\textrm{Ker}}\big[\emph{\textrm{Re }}F(\emph{\textrm{Im }}F)^j \big]\Big)\cap \rr^{2n}\\
=\big\{(x,\xi) \in \rr^{2n} : \forall 0 \leq j \leq \inf(k,n-1),\ QQ_{\infty}^{-1}B^jx=0,\ Q(B^T)^j\xi=0\big\}.
\end{multline*}
\end{lemma}

\medskip

\begin{proof}
We begin by proving by induction that for all $k \geq 0$,
\begin{multline}\label{pav23.1}
\Big(\bigcap_{j=0}^{k}\textrm{Ker}\big[\textrm{Re }F(\textrm{Im }F)^j \big]\Big)\cap \rr^{2n}\\
=\big\{(x,\xi) \in \rr^{2n} : \forall 0 \leq j \leq k,\ QQ_{\infty}^{-1}B^jx=0,\ Q(B^T)^j\xi=0\big\}.
\end{multline}
We notice from (\ref{pav23}) that the formula (\ref{pav23.1}) holds for $k=0$. If for $k \geq 0$,
\begin{multline*}
\Big(\bigcap_{j=0}^{k}\textrm{Ker}\big[\textrm{Re }F(\textrm{Im }F)^j \big]\Big)\cap \rr^{2n}\\
=\big\{(x,\xi) \in \rr^{2n} : \forall 0 \leq j \leq k,\ QQ_{\infty}^{-1}B^jx=0,\ Q(B^T)^j\xi=0\big\},
\end{multline*}
it follows from (\ref{pav21}) and (\ref{pav22}) that 
\begin{multline*}
\Big(\bigcap_{j=0}^{k+1}\textrm{Ker}\big[\textrm{Re }F(\textrm{Im }F)^j \big]\Big)\cap \rr^{2n}
=\big\{(x,\xi) \in \rr^{2n} : \forall 0 \leq j \leq k,\ QQ_{\infty}^{-1}B^jx=0,\\ Q(B^T)^j\xi=0, \ QQ_{\infty}^{-1}(QQ_{\infty}^{-1}+2B)^{k+1}x=0,\ Q(Q_{\infty}^{-1}Q+2B^T)^{k+1}\xi=0\big\}.
\end{multline*}
We deduce by expanding the products that 
$$QQ_{\infty}^{-1}(QQ_{\infty}^{-1}+2B)^{k+1}x=2^{k+1}QQ_{\infty}^{-1}B^{k+1}x,\quad Q(Q_{\infty}^{-1}Q+2B^T)^{k+1}\xi=2^{k+1}Q(B^T)^{k+1}\xi,$$
when
$$\forall 0 \leq j \leq k,\quad QQ_{\infty}^{-1}B^jx=0,\quad Q(B^T)^j\xi=0.$$
This proves that 
\begin{multline*}
\Big(\bigcap_{j=0}^{k+1}\textrm{Ker}\big[\textrm{Re }F(\textrm{Im }F)^j \big]\Big)\cap \rr^{2n}\\
=\big\{(x,\xi) \in \rr^{2n} : \forall 0 \leq j \leq k+1,\ QQ_{\infty}^{-1}B^jx=0,\ Q(B^T)^j\xi=0\big\}.
\end{multline*}
This shows that the formula (\ref{pav23.1}) holds for all $k \geq 0$. Next, we deduce from (\ref{pav23.1}) and the Cayley-Hamilton theorem 
$$\chi_B(B)=\chi_{B^T}(B^T)=0,$$
where $\chi_B$ and $\chi_{B^T}$ denote respectively the characteristic polynomials of the $n \times n$-matrices $B$ and $B^T$, that for all $k \geq 0$,
\begin{multline*}
\Big(\bigcap_{j=0}^{k}\textrm{Ker}\big[\textrm{Re }F(\textrm{Im }F)^j \big]\Big)\cap \rr^{2n}\\
=\big\{(x,\xi) \in \rr^{2n} : \forall 0 \leq j \leq \inf(k,n-1),\ QQ_{\infty}^{-1}B^jx=0,\ Q(B^T)^j\xi=0\big\},
\end{multline*}
since $\deg \chi_B=\deg \chi_{B^T}=n$. This ends the proof of Lemma~\ref{l1}.
\end{proof}

\medskip

We establish the following result:

\medskip

\begin{lemma}\label{l2}
For all $k \geq 0$, the following equivalence holds
$$\Big(\bigcap_{j=0}^{k}\emph{\textrm{Ker}}\big[\emph{\textrm{Re }}F(\emph{\textrm{Im }}F)^j \big]\Big)\cap \rr^{2n}=\{0\} \Longleftrightarrow \Big(\bigcap_{j=0}^{\inf(k,n-1)}\emph{\textrm{Ker}}\big[Q(B^T)^j\big]\Big)\cap \rr^{n}=\{0\}.$$
\end{lemma}

\medskip

\begin{proof}
According to Lemma~\ref{l1}, we have already proved that 
$$\Big(\bigcap_{j=0}^{k}\textrm{Ker}\big[\textrm{Re }F(\textrm{Im }F)^j \big]\Big)\cap \rr^{2n}=\{0\} \Longrightarrow \Big(\bigcap_{j=0}^{\inf(k,n-1)}\textrm{Ker}\big[Q(B^T)^j\big]\Big)\cap \rr^{n}=\{0\}.$$
Conversely, if 
$$\Big(\bigcap_{j=0}^{\inf(k,n-1)}\textrm{Ker}\big[Q(B^T)^j\big]\Big)\cap \rr^{n}=\{0\},$$
it follows from Lemma~\ref{l1} that it is sufficient to prove that 
$$\Big(\bigcap_{j=0}^{\inf(k,n-1)}\textrm{Ker}\big[QQ_{\infty}^{-1}B^j\big]\Big)\cap \rr^{n}=\Big(\bigcap_{j=0}^{\inf(k,n-1)}\textrm{Ker}\big[Q(Q_{\infty}^{-1}BQ_{\infty})^jQ_{\infty}^{-1}\big]\Big)\cap \rr^{n}=\{0\},$$
i.e.
$$\Big(\bigcap_{j=0}^{\inf(k,n-1)}\textrm{Ker}\big[Q(Q_{\infty}^{-1}BQ_{\infty})^j\big]\Big)\cap \rr^{n}=\{0\},$$
since $Q_{\infty}^{-1} \in M_n(\rr)$.
According to (\ref{pav10}), it is sufficient to show that  
$$\Big(\bigcap_{j=0}^{\inf(k,n-1)}\textrm{Ker}\big[Q(Q_{\infty}^{-1}Q+B^T)^j\big]\Big)\cap \rr^{n}=\{0\}.$$
On the other hand, we easily check by induction as in Lemma~\ref{l1} that for all $k \geq 0$,  
$$\Big(\bigcap_{j=0}^{k}\textrm{Ker}\big[Q(Q_{\infty}^{-1}Q+B^T)^j\big]\Big)\cap \rr^{n}=\Big(\bigcap_{j=0}^{k}\textrm{Ker}\big[Q(B^T)^j\big]\Big)\cap \rr^{n}.$$
This ends the proof of Lemma~\ref{l2}.
\end{proof}

We deduce from (\ref{h1bis}) and Lemma~\ref{l2} that the singular space of the quadratic form (\ref{pav17}) is zero if and only if the following intersection of kernels is zero
\begin{equation}\label{pav24}
S=\{0\} \Longleftrightarrow  \Big(\bigcap_{j=0}^{n-1}\textrm{Ker}\big[Q(B^T)^j\big]\Big)\cap \rr^{n}=\{0\} .
\end{equation}
On the other hand, the Kalman rank condition which holds since the Ornstein-Uhlenbeck operator (\ref{jen0}) is assumed to be hypoelliptic, reads as 
$$\textrm{Rank}[B|Q^{\frac{1}{2}}]=\textrm{Rank}[Q^{\frac{1}{2}},BQ^{\frac{1}{2}},\dots, B^{n-1}Q^{\frac{1}{2}}]=n,$$ 
where the $n\times n^2$ matrix $[Q^{\frac{1}{2}},BQ^{\frac{1}{2}},\dots, B^{n-1}Q^{\frac{1}{2}}]$ is obtained by writing consecutively the columns of the matrices $B^jQ^{\frac{1}{2}}$.
Writing $(C_1,C_2,...,C_n)$ the $n$ columns of the real matrix
$$\left[\begin{array}{c}
Q^{\frac{1}{2}}\\ Q^{\frac{1}{2}}B^T \\ \vdots\\ Q^{\frac{1}{2}}(B^T)^{n-1}
\end{array}\right],$$
the Kalman rank condition is equivalent to the linear independence of the column vectors $(C_1,C_2,...,C_n)$, which is also equivalent to the condition
$$\Big(\bigcap_{j=0}^{n-1}\textrm{Ker}\big[Q^{\frac{1}{2}}(B^T)^j\big]\Big)\cap \rr^{n}=\Big(\bigcap_{j=0}^{n-1}\textrm{Ker}\big[Q(B^T)^j\big]\Big)\cap \rr^{n}=\{0\}.$$
According to (\ref{pav24}), this shows that the hypoellipticity of the Ornstein-Uhlenbeck operator (\ref{jen0}) implies that the singular space of the quadratic operator (\ref{pav5}) is zero 
\begin{equation}\label{pav50}
S=\{0\}.
\end{equation} 
Regarding the subelliptic properties of the quadratic operator (\ref{pav5}), it is interesting to notice from Lemma~\ref{l2} that the smallest integer $0 \leq k_0 \leq 2n-1$ satisfying 
\begin{equation}\label{pav30}
\Big(\bigcap_{j=0}^{k_0}\textrm{Ker}
\big[\textrm{Re }F(\textrm{Im }F)^j \big]\Big)\cap \rr^{2n}=\{0\},
\end{equation}
corresponds exactly to the smallest integer $0 \leq k_0 \leq n-1$ satisfying
\begin{equation}\label{pav31}
\textrm{Rank}[Q^{\frac{1}{2}},BQ^{\frac{1}{2}},\dots, B^{k_0}Q^{\frac{1}{2}}]=n.
\end{equation}

\subsection{Condition for normality}
The following lemma gives a necessary condition for the operator $\mathscr{L}$ to be normal:

\medskip

\begin{lemma}\label{normal8}
Let 
$$P=\frac{1}{2}\emph{\textrm{Tr}}(Q\nabla_x^2)+\langle Bx,\nabla_x\rangle, \quad x \in \rr^n,$$ 
be a hypoelliptic Ornstein-Uhlenbeck operator, which admits the invariant measure $d\mu(x)=\rho(x) dx$.
Let $\mathscr{L}$ be the quadratic operator (\ref{pav4.5}) and $\mathscr{L}^*$ its $L^2(\rr^n,dx)$-adjoint.
Then, the commutator $[\mathscr{L},\mathscr{L}^*]$ is given by the quadratic operator
$$[\mathscr{L},\mathscr{L}^*]=\big\langle (QQ_{\infty}^{-1}Q+2QB^T)D_x,D_x \big\rangle+\frac{1}{4}\big\langle (QQ_{\infty}^{-1}Q+2QB^T)Q_{\infty}^{-1}x,Q_{\infty}^{-1}x\big\rangle,$$
with $D_x=i^{-1}\nabla_x$.
The operator $\mathscr{L}$ commutes with its $L^2$-adjoint $\mathscr{L}^*$,
$$\forall u \in \mathscr{S}(\rr^n), \quad [\mathscr{L},\mathscr{L}^*]u=0,$$
when acting on $\mathscr{S}(\rr^n)$, 
if and only if the following commutator is zero 
$$[Q_{\infty}^{-1}Q,B^T]=0.$$ 
Furthermore, if the operator $\mathscr{L}$ is normal, i.e.,
\begin{equation}\label{normal2}
\forall z \not\in \sigma(\mathscr{L}), \forall \zeta \not\in \sigma(\mathscr{L}^*), \quad (\mathscr{L}-z)^{-1}(\mathscr{L}^*-\zeta)^{-1}=(\mathscr{L}^*-\zeta)^{-1}(\mathscr{L}-z)^{-1},
\end{equation} 
then the condition $[Q_{\infty}^{-1}Q,B^T]=0$ holds. Moreover, the operator $\mathscr{L}$ can be normal only if the symmetric matrix $Q$ is positive definite.
\end{lemma}

\medskip
 
\begin{proof}
A direct computation or the use of Weyl calculus (see e.g. \cite{hormander}, Theorem~18.5.4) shows that the commutator $[\mathscr{L},\mathscr{L}^*]$ is a quadratic operator whose Weyl symbol is exactly given by the Poisson bracket
$$\frac{1}{i}\{q,\overline{q}\}=2\{\textrm{Im }q,\textrm{Re }q\}=2\langle \nabla_{\xi}\textrm{Im }q,\nabla_{x}\textrm{Re }q\rangle-2\langle \nabla_{x}\textrm{Im }q,\nabla_{\xi}\textrm{Re }q\rangle,$$
where $q$ denotes the Weyl symbol of $\mathscr{L}$. By using from (\ref{pav10}) that 
$$QQ_{\infty}^{-1}B=-QB^TQ_{\infty}^{-1}-QQ_{\infty}^{-1}QQ_{\infty}^{-1},$$
we deduce from (\ref{pav6}) that 
\begin{multline}\label{ql11}
2\{\textrm{Im }q,\textrm{Re }q\}=-\frac{1}{2}\Big\langle\Big(\frac{1}{2}QQ_{\infty}^{-1}+B\Big)x,Q_{\infty}^{-1}QQ_{\infty}^{-1}x\Big\rangle
+2\Big\langle\Big(\frac{1}{2}Q_{\infty}^{-1}Q+B^T\Big)\xi,Q\xi\Big\rangle\\
=\frac{1}{2}\Big\langle\Big(\frac{1}{2}QQ_{\infty}^{-1}Q+QB^T\Big)Q_{\infty}^{-1}x,Q_{\infty}^{-1}x\Big\rangle
+2\Big\langle\Big(\frac{1}{2}QQ_{\infty}^{-1}Q+QB^T\Big)\xi,\xi\Big\rangle,
\end{multline}
since $Q$ and $Q_{\infty}$ are symmetric. It follows that the commutator $[\mathscr{L},\mathscr{L}^*]$ is zero as an operator acting on the Schwartz space $\mathscr{S}(\rr^n)$,
\begin{equation}\label{normal1}
\forall u \in \mathscr{S}(\rr^n), \quad \mathscr{L}\mathscr{L}^*u=\mathscr{L}^*\mathscr{L}u,
\end{equation}
if and only if the following condition holds
$$QQ_{\infty}^{-1}Q+QB^T+BQ=0.$$
By using again from (\ref{pav10}) that
$$QQ_{\infty}^{-1}=-B-Q_{\infty}B^TQ_{\infty}^{-1},$$
this condition is equivalent to 
\begin{equation}\label{li1}
-Q_{\infty}B^TQ_{\infty}^{-1}Q+QB^T=0 \Longleftrightarrow [Q_{\infty}^{-1}Q,B^T]=0.
\end{equation}
By using the fact that the singular space is zero $S=\{0\}$, we deduce from (\ref{pav24}) and (\ref{li1}) that 
\begin{eqnarray*}
\{0\} &=& \Big(\bigcap_{j=0}^{n-1}\textrm{Ker}\big[Q(B^T)^j\big]\Big)\cap \rr^{n}=\Big(\bigcap_{j=0}^{n-1}\textrm{Ker}\big[Q_{\infty}^{-1}Q(B^T)^j\big]\Big)\cap \rr^{n}
\\ &=&
\Big(\bigcap_{j=0}^{n-1}\textrm{Ker}\big[(B^T)^jQ_{\infty}^{-1}Q\big]\Big)\cap \rr^{n}=\textrm{Ker}(Q) \cap \rr^{n},
\end{eqnarray*}
when the condition (\ref{li1}) holds.
It follows that the operator $\mathscr{L}$ satisfies the condition (\ref{normal1}) only if the real symmetric matrix $Q$ is positive definite. It remains to check the condition (\ref{normal1}) holds if the operator $\mathscr{L}$ is normal. We observe from (\ref{jk1}) and~\cite{kato} (Chap.~III, Sect.~6.6, Thm~6.22)  that $\sigma(\mathscr{L}) \subset \cc_+$ and $\sigma(\mathscr{L}^*) \subset \cc_+$, implying that $\sigma(\mathscr{L}) \neq \cc$ and $\sigma(\mathscr{L}^*) \neq \cc$. Let $z \not\in \sigma(\mathscr{L})$ and $\zeta \not\in \sigma(\mathscr{L}^*)$. We deduce from the inclusions 
$$\mathscr{S}(\rr^n) \subset D(\mathscr{L})=\{u \in L^2 : \mathscr{L}u \in L^2\},\quad \mathscr{S}(\rr^n) \subset D(\mathscr{L}^*)=\{u \in L^2 : \mathscr{L}^*u \in L^2\},$$ 
that 
$$\forall u \in \mathscr{S}(\rr^n), \exists v \in L^2, \ u=(\mathscr{L}-z)^{-1}v, \quad 
\forall u \in \mathscr{S}(\rr^n), \exists w \in L^2, \ u=(\mathscr{L}^*-\zeta)^{-1}w.$$
Since $(\mathscr{L}-z)u=v \in \mathscr{S}(\rr^n)$ and $(\mathscr{L}^*-\zeta)u=w \in \mathscr{S}(\rr^n)$, when $u \in \mathscr{S}(\rr^n)$, this implies that 
$$\forall u \in \mathscr{S}(\rr^n), \exists v \in \mathscr{S}(\rr^n), \ u=(\mathscr{L}-z)^{-1}v; \ \forall u \in \mathscr{S}(\rr^n), \exists w \in \mathscr{S}(\rr^n), \ u=(\mathscr{L}^*-\zeta)^{-1}w.$$
It follows that 
\begin{equation}\label{nh1}
\forall u \in \mathscr{S}(\rr^n), \exists v \in \mathscr{S}(\rr^n), \quad u=(\mathscr{L}-z)^{-1}(\mathscr{L}^*-\zeta)^{-1}v=(\mathscr{L}^*-\zeta)^{-1}(\mathscr{L}-z)^{-1}v,
\end{equation}
when the operator $\mathscr{L}$ is normal. We deduce from (\ref{nh1}) that for all $u \in \mathscr{S}(\rr^n)$,
$$(\mathscr{L}^*-\zeta)(\mathscr{L}-z)u=v=(\mathscr{L}-z)(\mathscr{L}^*-\zeta)u.$$
This implies that the operator $\mathscr{L}$ commutes with its $L^2$-adjoint $\mathscr{L}^*$,
$$\forall u \in \mathscr{S}(\rr^n), \quad [\mathscr{L},\mathscr{L}^*]u=0,$$
when acting on $\mathscr{S}(\rr^n)$. This ends the proof of Lemma~\ref{normal8}.
\end{proof}

\subsection{Computation of the eigenvalues of the Hamilton map}\label{res}
We aim at computing the eigenvalues of the Hamilton map (\ref{F}). After conjugating with the invertible matrix
$$\mathcal{G}=\lv\begin{array}{cc}
\frac{1}{\sqrt{2}}Q_{\infty}^{-\frac{1}{2}}& 0 \\
0 & \sqrt{2}Q_{\infty}^{\frac{1}{2}}       
\end{array}\rv,
$$
we notice from (\ref{F}) that  
\begin{equation}\label{pav32}
\mathcal{G}F\mathcal{G}^{-1}=\frac{1}{4}\lv\begin{array}{cc}
 -i\big(Q_{\infty}^{-\frac{1}{2}}QQ_{\infty}^{-\frac{1}{2}}+2Q_{\infty}^{-\frac{1}{2}}BQ_{\infty}^{\frac{1}{2}}\big)&Q_{\infty}^{-\frac{1}{2}}QQ_{\infty}^{-\frac{1}{2}} \\
 -Q_{\infty}^{-\frac{1}{2}}QQ_{\infty}^{-\frac{1}{2}}  & i\big(Q_{\infty}^{-\frac{1}{2}}QQ_{\infty}^{-\frac{1}{2}}+2Q_{\infty}^{-\frac{1}{2}}BQ_{\infty}^{\frac{1}{2}}\big)^T       
\end{array}\rv.
\end{equation}
Setting 
\begin{equation}\label{pav33}
M=-\frac{i}{2}Q_{\infty}^{-\frac{1}{2}}BQ_{\infty}^{\frac{1}{2}} \in M_n(\CC),
\end{equation}
we observe from (\ref{pav10}) that 
\begin{equation}\label{pav34}
M^T=-\frac{i}{2}Q_{\infty}^{\frac{1}{2}}B^TQ_{\infty}^{-\frac{1}{2}}=\frac{i}{2}(Q_{\infty}^{-\frac{1}{2}}QQ_{\infty}^{-\frac{1}{2}}+Q_{\infty}^{-\frac{1}{2}}BQ_{\infty}^{\frac{1}{2}}),
\end{equation}
since
$$Q_{\infty}^{-\frac{1}{2}}QQ_{\infty}^{-\frac{1}{2}}+Q_{\infty}^{-\frac{1}{2}}BQ_{\infty}^{\frac{1}{2}}+Q_{\infty}^{\frac{1}{2}}B^TQ_{\infty}^{-\frac{1}{2}}=0.$$
It follows from (\ref{pav32}), (\ref{pav33}) and (\ref{pav34}) that 
\begin{equation}\label{pav35}
\mathcal{G}F\mathcal{G}^{-1}=\mathcal{M},
\end{equation}
with
\begin{equation}\label{pav36}
\mathcal{M}=\lv\begin{array}{cc}
 \frac{M-M^T}{2}&-i\frac{M+M^T}{2} \\
 i\frac{M+M^T}{2}  & \frac{M-M^T}{2}   
\end{array}\rv.
\end{equation}
We use the following lemma:

\medskip

\begin{lemma}\label{pav37}
If 
$$\mathcal{M}=\lv\begin{array}{cc}
 \frac{M-M^T}{2}&-i\frac{M+M^T}{2} \\
 i\frac{M+M^T}{2}  & \frac{M-M^T}{2}   
\end{array}\rv \in M_{2n}(\CC),$$
with $M \in M_n(\CC)$, we have for all $X \in \CC^n$, $\mu \in \CC$ and $k \geq 0$,
$$(\mathcal{M}-\mu)^k\left[\begin{array}{c}
 X \\
 iX    
\end{array}\right]=\left[\begin{array}{c}
 (M-\mu)^kX \\
 i(M-\mu)^kX    
\end{array}\right], \ \ (\mathcal{M}-\mu)^k\left[\begin{array}{c}
 X \\
 -iX    
\end{array}\right]=(-1)^k\left[\begin{array}{c}
 (M^T+\mu)^kX \\
 -i(M^T+\mu)^kX    
\end{array}\right].$$
The spectrum of $\mathcal{M}$ is the union of the spectra of $M$ and $-M$,
$$\sigma(\mathcal{M})=\sigma(M) \cup \sigma(-M).$$
More precisely, the algebraic multiplicity of $\mu$ as an eigenvalue of $\mathcal{M}$ is equal to the sum of the algebraic multiplicity of $\mu$ as an eigenvalue of $M$ and the algebraic multiplicity of $\mu$ as an eigenvalue of $-M$.
\end{lemma}

\medskip

\begin{proof}
It is sufficient to check that 
$$(\mathcal{M}-\mu)\left[\begin{array}{c}
 X \\
 iX    
\end{array}\right]=\left[\begin{array}{c}
 (M-\mu)X \\
 i(M-\mu)X    
\end{array}\right], \ \ (\mathcal{M}-\mu)\left[\begin{array}{c}
 X \\
 -iX    
\end{array}\right]=-\left[\begin{array}{c}
 (M^T+\mu)X \\
 -i(M^T+\mu)X    
\end{array}\right],$$
and iterate these two formulas in order to prove the first assertion. Then, we use that 
\begin{equation}\label{pav38}
\CC^{2n}=F_1 \oplus F_2,
\end{equation}
where
\begin{equation}\label{pav39}
F_1=\{(X,iX) : X \in \CC^n\}, \quad F_2=\{(X,-iX) : X \in \CC^n\}.
\end{equation}
Let $\mathcal{B}=(e_1,...,e_n)$ be a basis of $\CC^n$ composed by generalized eigenvectors of the matrix $M \in M_n(\CC)$, and $\mathcal{C}=(\varepsilon_1,...,\varepsilon_n)$ be a basis of $\CC^n$ composed by generalized eigenvectors of the matrix $M^T \in M_n(\CC)$. It follows from the first assertion of this lemma that 
$$\tilde{\mathcal{B}}=(\tilde{e}_1,...,\tilde{e}_n),\ \textrm{ with }\ \tilde{e}_j=\left[\begin{array}{c}
 e_j \\
 ie_j    
\end{array}\right],$$
 is a basis of the vector subspace $F_1$ composed by generalized eigenvectors of the matrix $\mathcal{M} \in M_{2n}(\CC)$, and that 
 $$\tilde{\mathcal{C}}=(\tilde{\varepsilon}_1,...,\tilde{\varepsilon}_n),\ \textrm{ with }\ \tilde{\varepsilon}_j=\left[\begin{array}{c}
 \varepsilon_j \\
 -i\varepsilon_j    
\end{array}\right],$$
 is a basis of the vector subspace $F_2$ composed by generalized eigenvectors of the matrix $\mathcal{M} \in M_{2n}(\CC)$. $\tilde{\mathcal{B}} \cup \tilde{\mathcal{C}}$ is therefore a basis of $\CC^{2n}$ composed by generalized eigenvectors of the matrix $\mathcal{M} \in M_{2n}(\CC)$. We deduce that 
$$\sigma(\mathcal{M})=\sigma(M) \cup \sigma(-M^T)=\sigma(M) \cup \sigma(-M),$$
since the eigenvalues of $M$ and $M^T$ agree with algebraic multiplicities.
Furthermore, the algebraic multiplicity of $\mu$ as an eigenvalue of $\mathcal{M}$ is equal to the sum of the algebraic multiplicity of $\mu$ as an eigenvalue of $M$ and the algebraic multiplicity of $\mu$ as an eigenvalue of $-M$.
This ends the proof of Lemma~\ref{pav37}.
\end{proof}

We deduce from (\ref{pav33}), (\ref{pav35}), (\ref{pav36}) and Lemma~\ref{pav37} the following result:

\medskip

\begin{corollary}\label{coro1}
The spectrum of the Hamilton map 
$$F=\lv\begin{array}{cc}
 -\frac{i}{4}(QQ_{\infty}^{-1}+2B)&\frac{1}{2}Q \\
 -\frac{1}{8}Q_{\infty}^{-1}QQ_{\infty}^{-1}  & \frac{i}{4}(QQ_{\infty}^{-1}+2B)^T       
\end{array}\rv,$$
of the quadratic operator 
$$\mathscr{L}=\frac{1}{2}|Q^{\frac{1}{2}}D_x|^2+\frac{1}{8}|Q^{\frac{1}{2}}Q_{\infty}^{-1}x|^2-i\Big\langle\Big(\frac{1}{2}QQ_{\infty}^{-1}+B\Big)x,D_x\Big\rangle,$$
is given by
$$\sigma(F)=\sigma\Big(-\frac{i}{2}B\Big) \cup \sigma\Big(\frac{i}{2}B\Big).$$
More precisely, the algebraic multiplicity of $\mu$ as an eigenvalue of $F$ is equal to the sum of the algebraic multiplicity of $\mu$ as an eigenvalue of $\frac{i}{2}B$ and the algebraic multiplicity of $\mu$ as an eigenvalue of $-\frac{i}{2}B$.
\end{corollary}

\subsection{Numerical range of the quadratic operator $\mathscr{L}$}
The following lemma shows that the numerical range $\Sigma(q)=\overline{q(\rr^{2n})}$ of the Weyl symbol (\ref{pav6}) is the closed right half-plane of $\cc$, when the matrix $Q$ is degenerate:

\medskip

\begin{lemma}\label{jk1b}
When the symmetric matrix $Q$ is degenerate, then the numerical range $\Sigma(q)=\overline{q(\rr^{2n})}$ of the quadratic symbol
$$q(x,\xi)=\frac{1}{2}|Q^{\frac{1}{2}}\xi|^2+\frac{1}{8}|Q^{\frac{1}{2}}Q_{\infty}^{-1}x|^2-i\Big\langle\Big(\frac{1}{2}QQ_{\infty}^{-1}+B\Big)x,\xi\Big\rangle,$$
is equal to 
$$\Sigma(q)=\{z \in \cc :  \emph{\textrm{Re }}z \geq 0\}.$$
Furthermore, for all $z \in \cc_+=\{z \in \cc : \emph{\textrm{Re }}z>0\}$, there exists $(x_0,\xi_0) \in \rr^{2n}$ such that $z=q(x_0,\xi_0)$ and
$$\{\emph{\textrm{Re }}q,\emph{\textrm{Im }}q\}(x_0,\xi_0)=\langle \nabla_{\xi} \emph{\textrm{Re }}q,\nabla_{x} \emph{\textrm{Im }}q\rangle (x_0,\xi_0)-\langle \nabla_{x} \emph{\textrm{Re }}q,\nabla_{\xi} \emph{\textrm{Im }}q\rangle (x_0,\xi_0) <0,$$
where $\{\emph{\textrm{Re }}q,\emph{\textrm{Im }}q\}$ denotes the Poisson bracket of the symbols $\emph{\textrm{Re }}q$ and $\emph{\textrm{Im }}q$.
\end{lemma}

\medskip

\begin{proof}
Since $Q$ is a degenerate real symmetric matrix, we can find $\zeta_0 \in \rr^n$, $\zeta_0 \neq 0$ such that $\zeta_0 \in \textrm{Ker }Q$. We consider the non-negative analytic function
\begin{equation}\label{ql5}
\forall t \geq 0, \quad F_0(t)=\langle Qe^{tB^T}\zeta_0,e^{tB^T}\zeta_0 \rangle=|Q^{\frac{1}{2}}e^{tB^T}\zeta_0|^2 \geq 0.
\end{equation}
This function vanishes at zero $F_0(0)=0$, because $\zeta_0 \in \textrm{Ker }Q$. On the other hand, we notice that the function $F_0$ cannot be identically equal to zero since 
$$\forall t>0, \quad 0<|Q_t^{\frac{1}{2}}\zeta_0|^2=\int_0^tF_0(s)ds,$$
according to (\ref{pav0}). By analyticity, this implies that the function $F_0$ cannot be flat at zero. It follows from (\ref{ql5}) that there exists $j_0 \geq 1$ such that 
\begin{equation}\label{ql1}
F_0(t)=\frac{F_0^{(2j_0)}(0)}{(2j_0)!}t^{2j_0}+o(t^{2j_0}), \qquad F_0'(t)=\frac{F_0^{(2j_0)}(0)}{(2j_0-1)!}t^{2j_0-1}+o(t^{2j_0-1}),
\end{equation}
with $F_0^{(2j_0)}(0)>0$, when $t \to 0$. We deduce from (\ref{ql1}) that there exists $\delta>0$ such that
\begin{equation}\label{ql2}
\forall 0<t \leq \delta, \quad  F_0(t)>0, \qquad  \forall 0<t \leq \delta, \quad  F_0'(t)>0.
\end{equation}
Setting 
\begin{equation}\label{ql6}
\forall t \in \rr, \quad G(t)=\langle BQ_{\infty}\zeta_0,e^{tB^T}\zeta_0\rangle,
\end{equation}
we can find $0<t_0 \leq \delta$ such that $G(t_0) \neq 0$. Indeed, this directly follows from the continuity of $G$ if $G(0) \neq 0$. On the other hand, if $G(0)=0$, it follows from (\ref{pav1}) and (\ref{pav10}) that 
$$G'(0)=\langle BQ_{\infty}\zeta_0,B^T\zeta_0 \rangle=-\langle Q_{\infty}B^T\zeta_0,B^T\zeta_0 \rangle=-|Q_{\infty}^{\frac{1}{2}}B^T\zeta_0|^2<0,$$
because the matrices $Q_{\infty}$, $B^T$ are invertible and $\zeta_0 \in \textrm{Ker }Q$, $\zeta_0 \neq 0$. This implies that there exists $0<t_0\leq \delta$ such that $G(t_0)<0$. We can therefore choose $0<t_0 \leq \delta$ such that $F_0(t_0)>0$, $F_0'(t_0)>0$ and $G(t_0) \neq 0$. It follows from (\ref{ql5}) and (\ref{ql6}) that 
\begin{equation}\label{ql7}
F_0(t_0)=|Q^{\frac{1}{2}}e^{t_0B^T}\zeta_0|^2>0, \quad F_0'(t_0)=2\langle QB^Te^{t_0B^T}\zeta_0,e^{t_0B^T}\zeta_0 \rangle>0
\end{equation}
and
\begin{equation}\label{ql8}
G(t_0)=\langle BQ_{\infty}\zeta_0,e^{t_0B^T}\zeta_0\rangle \neq 0.
\end{equation}
Let $z \in \cc_+$. Setting
\begin{equation}\label{ql9}
\lambda=-\frac{\textrm{Im }z}{G(t_0)}\sqrt{\frac{F_0(t_0)}{2\textrm{Re }z}}, \qquad \mu=\sqrt{\frac{2\textrm{Re }z}{F_0(t_0)}}>0,
\end{equation}
we consider
\begin{equation}\label{ql10}
(x_0,\xi_0)=(\lambda Q_{\infty}\zeta_0,\mu e^{t_0B^T}\zeta_0) \in \rr^{2n}.
\end{equation}
It follows from (\ref{pav6}), (\ref{ql7}), (\ref{ql8}), (\ref{ql9}) and (\ref{ql10}) that 
\begin{multline}\label{ql11v}
q(x_0,\xi_0)=\frac{\mu^2}{2}|Q^{\frac{1}{2}}e^{t_0B^T}\zeta_0|^2+\frac{\lambda^2}{8}|Q^{\frac{1}{2}}\zeta_0|^2-i\lambda \mu \Big\langle\Big(\frac{1}{2}Q+BQ_{\infty}\Big)\zeta_0,e^{t_0B^T}\zeta_0\Big\rangle\\
=\frac{\mu^2}{2}|Q^{\frac{1}{2}}e^{t_0B^T}\zeta_0|^2-i\lambda \mu \langle BQ_{\infty}\zeta_0,e^{t_0B^T}\zeta_0\rangle
=\frac{\mu^2}{2}F_0(t_0)-i\lambda \mu G(t_0)=z,
\end{multline}
since $\zeta_0 \in \textrm{Ker }Q$. It implies that $\cc_+ \subset q(\rr^{2n})$. On the other hand, we have 
$$q(\rr^{2n}) \subset \{z \in \cc : \textrm{Re }z \geq 0\},$$
since the real part of $q$ is non-negative. It follows that 
$$\Sigma(q)=\overline{q(\rr^{2n})}=\{z \in \cc : \textrm{Re }z \geq 0\}.$$
We deduce from (\ref{ql11}) and (\ref{ql10}) that 
\begin{multline*}
\{\textrm{Re }q,\textrm{Im }q\}(x_0,\xi_0)=-\frac{1}{4}\Big\langle\Big(\frac{1}{2}QQ_{\infty}^{-1}Q+QB^T\Big)Q_{\infty}^{-1}x_0,Q_{\infty}^{-1}x_0\Big\rangle\\
-\Big\langle\Big(\frac{1}{2}QQ_{\infty}^{-1}Q+QB^T\Big)\xi_0,\xi_0\Big\rangle=-\mu^2\Big\langle\Big(\frac{1}{2}QQ_{\infty}^{-1}Q+QB^T\Big)e^{t_0B^T}\zeta_0,e^{t_0B^T}\zeta_0\Big\rangle,
\end{multline*}
because $Q$ is symmetric and $Q\zeta_0=0$. It follows from (\ref{ql7}) and (\ref{ql9}) that 
$$\{\textrm{Re }q,\textrm{Im }q\}(x_0,\xi_0)
=-\frac{\mu^2}{2} |Q_{\infty}^{-\frac{1}{2}}Qe^{t_0B^T}\zeta_0|^2-\frac{\mu^2}{2}F_0'(t_0)<0.$$ 
This ends the proof of Lemma~\ref{jk1b}.
\end{proof}

\bigskip
\noindent
{\bf Acknowledgements.}
The research of the second author is supported by the Engineering and Physical Sciences Research Council of the UK through Grants No. EP/J009636/1, EP/L024926/1, EP/L020564/1 and EP/L025159/1. The research of the third author is supported by the ANR NOSEVOL (Project: ANR 2011 BS01019 01). The figures in the article have been created using the software Fig4Tex developed by Lafranche and Martin~\cite{lafranche}.

\end{document}